\documentclass[12pt]{article}
\usepackage{amsmath,amsfonts,amssymb,amsthm}
\input amssym.def
\begin{document}
\def \Z{\Bbb Z}
\def \C{\Bbb C}
\def \R{\Bbb R}
\def \Q{\Bbb Q}
\def \N{\Bbb N}
\def \bR{\bf R}
\def \D{{\cal{D}}}
\def \E{{\cal{E}}}
\def \S{{\cal{S}}}
\def \R{{\cal{R}}}
\def \Y{{\cal{Y}}}
\def \wt{{\rm wt}}
\def \tr{{\rm tr}}
\def \span{{\rm span}}
\def \Res{{\rm Res}}
\def \End{{\rm End}\;}
\def \Ind {{\rm Ind}}
\def \Irr {{\rm Irr}}
\def \Aut{{\rm Aut}}
\def \Hom{{\rm Hom}}
\def \mod{{\rm mod}}
\def \ann{{\rm Ann}}
\def \ad{{\rm ad}}
\def \rank{{\rm rank}\;}

\def \<{\langle}
\def \>{\rangle}
\def \t{\tau }
\def \a{\alpha }
\def \e{{\bf e}}
\def \f{{\bf f}}
\def \l{\lambda }
\def \L{\Lambda }
\def \g{{\frak{g}}}
\def \h{{\hbar}}
\def \k{{\frak{k}}}
\def \sl{{\frak{sl}}}
\def \be{\begin{equation}\label}
\def \ee{\end{equation}}
\def \bex{\begin{example}\label}
\def \eex{\end{example}}
\def \bl{\begin{lem}\label}
\def \el{\end{lem}}
\def \bt{\begin{thm}\label}
\def \et{\end{thm}}
\def \bp{\begin{prop}\label}
\def \ep{\end{prop}}
\def \br{\begin{rem}\label}
\def \er{\end{rem}}
\def \bc{\begin{coro}\label}
\def \ec{\end{coro}}
\def \bd{\begin{de}\label}
\def \ed{\end{de}}
\def \pf{{\bf Proof. }}

\newcommand{\n}{\:^{\times}_{\times}\:}
\newcommand{\nno}{\nonumber}
\newcommand{\nord}{\mbox{\scriptsize ${\circ\atop\circ}$}}
\newtheorem{thm}{Theorem}[section]
\newtheorem{prop}[thm]{Proposition}
\newtheorem{coro}[thm]{Corollary}
\newtheorem{conj}[thm]{Conjecture}
\newtheorem{example}[thm]{Example}
\newtheorem{lem}[thm]{Lemma}
\newtheorem{rem}[thm]{Remark}
\newtheorem{de}[thm]{Definition}
\newtheorem{hy}[thm]{Hypothesis}
\makeatletter
\@addtoreset{equation}{section}
\def\theequation{\thesection.\arabic{equation}}
\makeatother
\makeatletter

\begin{center}{\Large \bf Modules-at-infinity for quantum vertex algebras}
\end{center}

\begin{center}
{Haisheng Li\footnote{Partially supported by NSA grant H98230-05-1-0018
and NSF grant DMS-0600189}\\
Department of Mathematical Sciences, Rutgers University, Camden, NJ
08102}
\end{center}

\begin{abstract}
This is a sequel to \cite{li-qva1} and \cite{li-qva2} in a series to
study vertex algebra-like structures arising from various algebras
such as quantum affine algebras and Yangians.  In this paper, we
study two versions of the double Yangian $DY_{\hbar}(sl_{2})$,
denoted by $DY_{q}(sl_{2})$ and $DY_{q}^{\infty}(sl_{2})$ with $q$ a
nonzero complex number.  For each nonzero complex number $q$, we
construct a quantum vertex algebra $V_{q}$ and prove that every
$DY_{q}(sl_{2})$-module is naturally a $V_{q}$-module. We also show
that $DY_{q}^{\infty}(sl_{2})$-modules are what we call
$V_{q}$-modules-at-infinity.  To achieve this goal, we study what we
call $\S$-local subsets and quasi-local subsets of $\Hom
(W,W((x^{-1})))$ for any vector space $W$, and we prove that any
$\S$-local subset generates a (weak) quantum vertex algebra and that
any quasi-local subset generates a vertex algebra with $W$ as a
(left) quasi module-at-infinity.  Using this result we associate the
Lie algebra of pseudo-differential operators on the circle with
vertex algebras in terms of quasi modules-at-infinity.
\end{abstract}

\section{Introduction}
This is a sequel to \cite{li-qva1} and \cite{li-qva2} in a series to
study vertex algebra-like structures arising from various algebras
such as quantum affine algebras and Yangians.  In \cite{li-qva1} and
\cite{li-qva2}, partially motivated by Etingof-Kazhdan's notion of
quantum vertex operator algebra over $\C[[\hbar]]$ (see \cite{ek}),
we formulated and studied a notion of quantum vertex algebra over
$\C$ and we established general constructions of (weak) quantum
vertex algebras and modules.  The general constructions were
illustrated by examples in which quantum vertex algebras were
constructed from certain Zamolodchikov-Faddeev-type algebras.  The
main goal of this paper is to establish a natural connection of
(centerless) double Yangians with quantum vertex algebras over $\C$.

For each finite-dimensional simple Lie algebra $\g$, Drinfeld 
introduced a Hopf algebra $Y(\g)$, called Yangian (see \cite{dr1}),
as a deformation of the universal enveloping algebra of the
Lie algebra $\g\otimes \C[t]$. The Yangian
double $DY_{\hbar}(\g)$ is a deformation of the algebra 
$U(\g\otimes \C[t,t^{-1}])$.
In the simplest case with $\g=sl_{2}$ (cf. \cite{sm}, \cite{kh}), the
following is one of the defining relations in terms of generating
functions
\begin{eqnarray}\label{e-introd-ee}
e(x_{1})e(x_{2})=\frac{x_{1}-x_{2}+\hbar}{x_{1}-x_{2}-\hbar}e(x_{2})e(x_{1}).
\end{eqnarray}
In this paper we study two versions of $DY_{\hbar}(sl_{2})$
with the formal parameter $\hbar$ being evaluated at 
a nonzero complex number $q$.
With a direct substitution, the defining relation (\ref{e-introd-ee}) becomes
\begin{eqnarray}\label{e-introd-qee}
e(x_{1})e(x_{2})=\frac{x_{1}-x_{2}+q}{x_{1}-x_{2}-q}e(x_{2})e(x_{1}),
\end{eqnarray}
where
\begin{eqnarray*}
\frac{x_{1}-x_{2}+q}{x_{1}-x_{2}-q}
=(x_{1}-x_{2}+q)\sum_{i\ge 0} q^{i}(x_{1}-x_{2})^{-i-1}
\in \C(((x_{1}-x_{2})^{-1}))\subset \C[x_{2}]((x_{1}^{-1})).
\end{eqnarray*}
In this way we get a version of $DY_{\hbar}(sl_{2})$, 
which we denote by $DY_{q}^{\infty}(sl_{2})$.
Notice that for a module $W$ of highest weight type where
the generating functions are elements of $\Hom (W,W((x)))$,
the expression
$$\frac{x_{1}-x_{2}+q}{x_{1}-x_{2}-q}e(x_{2})e(x_{1})$$
does not exist in general.
Because of this, $DY_{q}^{\infty}(sl_{2})$
admits only modules $W$ of lowest weight type 
where the generating functions
such as $e(x)$ are elements of $\Hom (W,W((x^{-1})))$.

Note that (quantum) vertex algebras and their modules are
modules of highest weight type in nature and
so far we have used only modules of highest weight type
 for various algebras
to construct (quantum) vertex algebras and their modules. 
Motivated by this,
we then consider another version of $DY_{\hbar}(sl_{2})$,
which we denote by $DY_{q}(sl_{2})$,
by expanding the same rational function
$\frac{x_{1}-x_{2}+q}{x_{1}-x_{2}-q}$
as follows:
$$\frac{x_{1}-x_{2}+q}{-q+x_{1}-x_{2}}
=-(x_{1}-x_{2}+q)\sum_{i\ge 0} q^{-i-1}(x_{1}-x_{2})^{i}
\in \C[[(x_{1}-x_{2})]]\subset \C[[x_{1},x_{2}]].$$
Contrary to the situation with $DY_{q}^{\infty}(sl_{2})$,
the algebra $DY_{q}(sl_{2})$
admits only modules of highest weight type, including vacuum modules.
For every nonzero complex number $q$, we construct a universal
vacuum $DY_{q}(sl_{2})$-module $V_{q}$ and
by applying the general construction theorems of
\cite{li-qva1} and \cite{li-qva2} we show that there exists
a canonical quantum vertex algebra structure on $V_{q}$
and that on any $DY_{q}(sl_{2})$-module there exists a canonical
$V_{q}$-module structure.

As it was mentioned before, the algebra $DY_{q}^{\infty}(sl_{2})$
admits only lowest-weight-type modules. This is also the case
for the Lie algebra
of pseudo-differential operators. Nevertheless, we hope to associate such
algebras with (quantum) vertex algebras in some natural way. Having
this in mind, we systematically study how to construct quantum
vertex algebras from suitable subsets of the space
$$\E^{o}(W)=\Hom (W,W((x^{-1})))$$
for a general vector space $W$, developing a general
theory analogous to that of \cite{li-qva1}.
For a vector space $W$,
a subset $T$ of $\E^{o}(W)$ is said to be $\S$-local if for any
$a(x),b(x)\in T$, there exist
$$a_{i}(x), b_{i}(x)\in T,\; f_{i}(x)\in \C(x),\; i=1,\dots,r$$
and a nonnegative integer $k$ such that
$$(x_{1}-x_{2})^{k}a(x_{1})b(x_{2})
=(x_{1}-x_{2})^{k}\sum_{i=1}^{r}\iota_{x,\infty}(f_{i})(x_{1}-x_{2})
b_{i}(x_{2})a_{i}(x_{1}),$$ where $\iota_{x,\infty}f_{i}(x)$ denotes
the formal Laurent series expansion of $f_{i}(x)$ at infinity. A
subset $T$ of $\E^{o}(W)$ is said to be quasi-local if for any
$a(x),b(x)\in T$, there exists a nonzero polynomial $p(x_{1},x_{2})$
such that
$$p(x_{1},x_{2})a(x_{1})b(x_{2})=p(x_{1},x_{2})b(x_{2})a(x_{1}).$$
We prove that any $\S$-local subset generates a weak quantum vertex
algebra and that any quasi-local subset of $\E^{o}(W)$ generates a
vertex algebra in a certain natural way. To describe the structure
on $W$ we formulate a notion of (left) quasi module-at-infinity for
a vertex algebra and for a weak quantum vertex algebra. For a vertex
algebra $V$, a (left) quasi $V$-module-at-infinity is a vector space
$W$ equipped with a linear map $Y_{W}$ from $V$ to $\Hom
(W,W((x^{-1})))$, satisfying the condition that $Y_{W}({\bf
1},x)=1_{W}$ (the identity operator on $W$) and that for $u,v\in V$,
there exists a nonzero polynomial $p(x_{1},x_{2})$ such that
\begin{eqnarray*}
& &p(x_{1},x_{2})Y_{W}(u,x_{1})Y_{W}(v,x_{2})
=p(x_{1},x_{2})Y_{W}(v,x_{2})Y_{W}(u,x_{1}),\\
& &p(x_{2}+x_{0},x_{2})Y_{W}(Y(u,x_{0})v,x_{2})
=\left(p(x_{1},x_{2})Y_{W}(u,x_{1})Y_{W}(v,x_{2})\right)|_{x_{1}=x_{2}+x_{0}}.
\end{eqnarray*}
For a module-at-infinity, the following opposite Jacobi identity holds for
$u,v\in V$:
\begin{eqnarray*}
& &x_{0}^{-1}\delta\left(\frac{x_{1}-x_{2}}{x_{0}}\right)
Y_{W}(v,x_{2})Y_{W}(u,x_{1})
-x_{0}^{-1}\delta\left(\frac{x_{2}-x_{1}}{-x_{0}}\right)
Y_{W}(u,x_{1})Y_{W}(v,x_{2})\\
& &\ \ \ \ =x_{2}^{-1}\delta\left(\frac{x_{1}-x_{0}}{x_{2}}\right)
Y_{W}(Y(u,x_{0})v,x_{2}).
\end{eqnarray*}
This notion of a left module-at-infinity for a vertex algebra $V$
coincides with the notion of a right module, which was
suggested in \cite{hl}.

As an application, we show that $DY_{q}(sl_{2})$-modules are canonical
modules-at-infinity for the quantum vertex algebra $V_{q}$.
We also show that lowest-weight type modules for
the Lie algebras of pseudo-differential operators on the circle 
are quasi modules-at-infinity for some
vertex algebras associated with the affine Lie
algebra of a certain infinite-dimensional Lie algebra.
In a sequel, we shall study general double Yangians $DY_{\hbar}(\g)$ and
their central extensions
$\widehat{DY_{\hbar}(\g)}$ (see \cite{kh}) 
in terms of quantum vertex algebras.

 This paper is organized as follows: In Section 2, we introduce a
 version of the double Yangian $DY_{\hbar}(sl_{2})$ and we associate
 it with quantum vertex algebras and modules.
  In Section 3, we study
 quasi-compatible subsets and prove that any quasi-compatible
 subset canonically generates a nonlocal vertex algebra.
  In Section 4, we study $\S$-local
 subsets and modules-at-infinity for quantum vertex algebras.
In Section 5, we study quasi local subsets and quasi
modules-at-infinity for vertex algebras.

\section{Associative algebra $DY_{q}(sl_{2})$ and quantum vertex algebras}

In this section,  we first recall the notions of weak quantum vertex
algebra and quantum vertex algebra and we then define an associative
algebra $DY_{q}(sl_{2})$ over $\C$ with $q$ an arbitrary nonzero
complex number, which is a version of the (centerless) double
Yangian $DY_{\hbar}(sl_{2})$, and we associate a quantum vertex
algebra to the algebra $DY_{q}(sl_{2})$.

We begin with the notion of nonlocal vertex algebra 
(\cite{kacbook2}, \cite{bk},
\cite{li-g1}). A {\em nonlocal vertex algebra} is a vector space
$V$, equipped with a linear map
\begin{eqnarray}
Y: V &\rightarrow & \Hom (V,V((x)))\subset (\End V)[[x,x^{-1}]]\nonumber\\
v&\mapsto& Y(v,x)=\sum_{n\in \Z}v_{n}x^{-n-1}\;\;\; (v_{n}\in \End V)
\end{eqnarray}
and equipped with a distinguished vector ${\bf 1}$, such that for
$v\in V$
\begin{eqnarray}
& &Y({\bf 1},x)v=v,\\
& &Y(v,x){\bf 1}\in V[[x]]\;\;\mbox{ and }\;\; \lim_{x\rightarrow
0}Y(v,x){\bf 1}=v
\end{eqnarray}
and such that for $u,v,w\in V$, there exists a nonnegative integer
$l$ such that
\begin{eqnarray}
(x_{0}+x_{2})^{l}Y(u,x_{0}+x_{2})Y(v,x_{2})w=
(x_{0}+x_{2})^{l}Y(Y(u,x_{0})v,x_{2})w.
\end{eqnarray}
For a nonlocal vertex algebra $V$, we have
\begin{eqnarray}\label{edproperty}
[{\cal{D}},Y(v,x)]=Y({\cal{D}}v,x)={d\over dx}Y(v,x) \;\;\;\mbox{
for }v\in V,
\end{eqnarray}
where  $\D$ is the linear operator on $V$, defined by
\begin{eqnarray*}
{\cal{D}}(v)=\left({d\over dx}Y(v,x){\bf 1}\right)|_{x=0}
\left(=v_{-2}{\bf 1}\right)
\;\;\;\mbox{ for }v\in V.
\end{eqnarray*}
Furthermore, for $v\in V$,
\begin{eqnarray}
& &e^{x{\cal{D}}}Y(v,x_{1})e^{-x{\cal{D}}}=Y(e^{xD}v,x_{1})=Y(v,x_{1}+x),
\label{econjugationformula1}\\
& &Y(v,x){\bf 1}=e^{x{\cal{D}}}v.\label{ecreationwithd}
\end{eqnarray}

A {\em weak quantum vertex algebra} (see \cite{li-qva1},
\cite{li-qva2}) is a vector space $V$ (over $\C$) equipped with a
distinguished vector ${\bf 1}$ and a linear map
\begin{eqnarray*}
Y:& &V\rightarrow \Hom (V,V((x)))\subset (\End V)[[x,x^{-1}]],\\
& &v\mapsto Y(v,x)=\sum_{n\in \Z}v_{n}x^{-n-1}\;\;
(\mbox{where }v_{n}\in \End V)
\end{eqnarray*}
 satisfying the condition that
\begin{eqnarray}
& &Y({\bf 1},x)v=v,\\
& &Y(v,x){\bf 1}\in V[[x]]\;\;\mbox{ and }\;\;\lim_{x\rightarrow
0}Y(v,x){\bf 1}=v\ \ \ \mbox{ for }v\in V,
\end{eqnarray}
 and that for any $u,v\in V$, there exist
$u^{(i)},v^{(i)}\in V,\; f_{i}(x)\in \C((x)),\; i=1,\dots,r$, such
that
\begin{eqnarray}\label{eS-jacobi}
&&x_{0}^{-1}\delta\left(\frac{x_{1}-x_{2}}{x_{0}}\right)Y(u,x_{1})Y(v,x_{2})
-x_{0}^{-1}\delta\left(\frac{x_{2}-x_{1}}{-x_{0}}\right)
\sum_{i=1}^{r}f_{i}(-x_{0})Y(v^{(i)},x_{2})Y(u^{(i)},x_{1})\nonumber\\
& &\hspace{2cm} = x_{2}^{-1}\delta\left(\frac{x_{1}-x_{0}}{x_{2}}\right)
Y(Y(u,x_{0})v,x_{2}).
\end{eqnarray}

In terms of the notion of nonlocal vertex algebra, a weak quantum
vertex algebra is simply a nonlocal vertex algebra that satisfies
the $\S$-locality (cf. \cite{ek}) in the sense that
for any $u,v\in V$, there exist
$u^{(i)},v^{(i)}\in V,\; f_{i}(x)\in \C((x)),\; i=1,\dots,r$, such
that
\begin{eqnarray}
(x_{1}-x_{2})^{k}Y(u,x_{1})Y(v,x_{2})
=(x_{1}-x_{2})^{k}\sum_{i=1}^{r}f_{i}(x_{2}-x_{1})
Y(v^{(i)},x_{2})Y(u^{(i)},x_{1})
\end{eqnarray}
for some nonnegative integer $k$ depending on $u$ and $v$.

The notion of quantum vertex algebra involves the notion of unitary rational
quantum Yang-Baxter operator, which here we recall:
A {\em unitary rational quantum Yang-Baxter operator}
on a vector space $H$ is a linear map
$\S(x): H\otimes H\rightarrow H\otimes H\otimes \C((x))$ such that
\begin{eqnarray*}
& &\S_{21}(-x)\S(x)=1,\\
& &\S_{12}(x)\S_{13}(x+z)\S_{23}(z)=\S_{23}(z)\S_{13}(x+z)\S_{12}(x).
\end{eqnarray*}
In this definition,  $\S_{21}(x)=\sigma_{12}\S(x)\sigma_{12}$, where
$\sigma_{12}$ is the flip map on $H\otimes H$ $(u\otimes v\mapsto
v\otimes u)$,
$$\S_{12}(x)=\S(x)\otimes 1: H\otimes H\otimes H
\rightarrow H\otimes H\otimes H\otimes \C((x)),$$
and $\S_{13}(x),\S_{23}(x)$ are defined accordingly.

A {\em quantum vertex algebra} (\cite{li-qva1}, cf. \cite{ek}) is a
weak quantum vertex algebra $V$ equipped with a unitary rational
quantum Yang-Baxter operator $\S(x): V\otimes V\rightarrow V\otimes
V\otimes \C((x))$ such that for $u,v\in V$, (\ref{eS-jacobi}) holds
with $\S(x)(v\otimes u)=\sum_{i=1}^{r}v^{(i)}\otimes u^{(i)}\otimes
f_{i}(x)$.

\br{rold-results}
{\em Recall from \cite{ek} that a nonlocal vertex algebra $V$
is {\em nondegenerate} if for every positive integer $n$,
the linear map
$$Z_{n}: V^{\otimes n}\otimes \C((x_{1}))\cdots ((x_{n}))
\rightarrow V((x_{1}))\cdots ((x_{n}))$$
defined by
$$Z_{n}(v_{(1)}\otimes \cdots \otimes v_{(n)}\otimes f)
=fY(v_{(1)},x_{1})\cdots Y(v_{(n)},x_{n}){\bf 1}$$ is injective. It
was proved (\cite{li-qva1}, Theorem 4.8, cf. \cite{ek}, Proposition
1.11) that if $V$ is a nondegenerate weak quantum vertex algebra,
then the $\S$-locality of vertex operators $Y(v,x)$ for $v\in V$
uniquely defines a linear map
$$\S(x): V\otimes V\rightarrow V\otimes V\otimes \C((x))$$
such that $V$ equipped with $\S(x)$ is a quantum vertex algebra and
$\S(x)$ is the unique rational quantum Yang-Baxter operator making
$V$ a quantum vertex algebra. In view of this, we shall use the term
``a nondegenerate quantum vertex algebra'' for a nondegenerate weak
quantum vertex algebra which is a quantum vertex algebra with the
canonical rational quantum Yang-Baxter operator.} \er

Let $V$ be a nonlocal vertex algebra. A {\em $V$-module} is a vector
space $W$ equipped with a linear map $Y_{W}: V\rightarrow \Hom
(W,W((x)))$ satisfying the condition that $Y_{W}({\bf 1},x)=1$
(the identity operator on
$W$) and for any $u,v\in V,\; w\in W$, there exists a nonnegative
integer $l$ such that
\begin{eqnarray}
(x_{0}+x_{2})^{l}Y_{W}(u,x_{0}+x_{2})Y_{W}(v,x_{2})w
=(x_{0}+x_{2})^{l}Y_{W}(Y(u,x_{0})v,x_{2})w.
\end{eqnarray}
Now, assume that $V$ is a weak quantum vertex algebra and let
$(W,Y_{W})$ be a module for $V$ viewed as a nonlocal vertex algebra.
It was proved (\cite{li-qva1}, Lemma 5.7) that for any $u,v\in V$,
\begin{eqnarray}
&&x_{0}^{-1}\delta\left(\frac{x_{1}-x_{2}}{x_{0}}\right)
Y_{W}(u,x_{1})Y_{W}(v,x_{2})\nonumber\\
& &\hspace{1cm}
-x_{0}^{-1}\delta\left(\frac{x_{2}-x_{1}}{-x_{0}}\right)
\sum_{i=1}^{r}f_{i}(-x_{0})Y_{W}(v^{(i)},x_{2})Y_{W}(u^{(i)},x_{1})\nonumber\\
&=&x_{2}^{-1}\delta\left(\frac{x_{1}-x_{0}}{x_{2}}\right)
Y_{W}(Y(u,x_{0})v,x_{2}),
\end{eqnarray}
where $u^{(i)},v^{(i)}\in V,\; f_{i}(x)\in \C((x))$ are the same as
those in (\ref{eS-jacobi}).

\br{rtwo-parts}
{\em Here, we recall a general construction of weak quantum
 vertex algebras from \cite{li-qva1}. Let
$W$ be any vector space (over $\C$) and set $\E(W)=\Hom (W,W((x)))$. A
subset $T$ of $\E(W)$ is $\S$-local if for any $a(x),b(x)\in T$,
there exist
$$f_{i}(x)\in \C((x)),\; a^{(i)}(x),b^{(i)}(x)\in T
\ \ (i=1,\dots,r)$$
(finitely many) such that
$$(x_{1}-x_{2})^{k}a(x_{1})b(x_{2})
=\sum_{i=1}^{r}(x_{1}-x_{2})^{k}f_{i}(x_{2}-x_{1})b^{(i)}(x_{2})a^{(i)}(x_{1})$$
for some nonnegative integer $k$. In this case, for any $n\in \Z$,
we define (see \cite{li-qva2})
$$a(x)_{n}b(x)=\Res_{x_{1}}\left((x_{1}-x)^{n}a(x_{1})b(x)
-(-x+x_{1})^{n}\sum_{i=1}^{r}f_{i}(x-x_{1})b^{(i)}(x)a^{(i)}(x_{1})\right).$$
(Note that it is well defined.) It was proved (\cite{li-qva1},
Theorem 5.8) that every $\S$-local subset of $\E(W)$ generates a
weak quantum vertex algebra with $W$ as a module with
$Y_{W}(a(x),x_{0})=a(x_{0})$. } \er

Next, we introduce a version of the double Yangian $DY_{\hbar}(sl_{2})$
(associated to the three-dimensional simple Lie algebra $sl_{2}$).
Let $T(sl_{2}\otimes \C[t,t^{-1}])$ denote the tensor algebra
over the vector space $sl_{2}\otimes \C[t,t^{-1}]$.
{}{\em From now on we shall simply use $T$ for this algebra.}
We equip $T$ with the $\Z$-grading which is uniquely defined by
$$\deg (u\otimes t^{n})=n\ \ \
\mbox{ for }u\in sl_{2},\; n\in \Z,$$
making $T$ a $\Z$-graded algebra $T=\coprod_{n\in \Z}T_{n}$.
For $n\in \Z$,
set
$$I[n]=\coprod_{m\ge n} T_{m}\subset T.$$
This defines a decreasing filtration
of $T=\cup_{n\in \Z}I[n]$ with $\cap_{n\in \Z} I[n]=0$.
Denote by $\overline{T}$ the completion of $T$ associated with this
filtration.

For $u\in sl_{2}$, set
$$u(x)=\sum_{n\in \Z}u(n)x^{-n-1},$$
where $u(n)=u\otimes t^{n}$. We also write
$sl_{2}(n)=sl_{2}\otimes t^{n}$ for $n\in \Z$.
Let $e,f,h$ be the standard Chevalley generators of $sl_{2}$.

\bd{dyangian-double-algebra}
{\em Let $q$ be a nonzero complex number. We define $DY_{q}(sl_{2})$
to be the quotient algebra of $\overline{T}$ modulo the following relations:
\begin{eqnarray*}
& &e(x_{1})e(x_{2})
=\frac{q+x_{1}-x_{2}}{-q+x_{1}-x_{2}}
e(x_{2})e(x_{1}),\\
& &f(x_{1})f(x_{2})
=\frac{-q+x_{1}-x_{2}}{q+x_{1}-x_{2}}
f(x_{2})f(x_{1}),\\
& &[e(x_{1}),f(x_{2})]
=x_{1}^{-1}\delta\left(\frac{x_{2}}{x_{1}}\right)h(x_{2}),\\
& &h(x_{1})e(x_{2})
=\frac{q+x_{1}-x_{2}}{-q+x_{1}-x_{2}}
e(x_{2})h(x_{1}),\\
& &h(x_{1})h(x_{2})=h(x_{2})h(x_{1}),\\
& &h(x_{1})f(x_{2})
=\frac{-q+x_{1}-x_{2}}{q+x_{1}-x_{2}}
f(x_{2})h(x_{1}),
\end{eqnarray*}
where it is understood that
$$(\pm q+x_{1}-x_{2})^{-1}=\sum_{i\in \N}(\pm q)^{-i-1}(x_{2}-x_{1})^{i}
\in \C[[x_{1},x_{2}]].$$}
\ed

It is straightforward to see that $DY_{q}(sl_{2})$ admits a (unique)
derivation $d$ such that
\begin{eqnarray}
[d,u(x)]=\frac{d}{dx}u(x)\ \ \ \mbox{ for }u\in sl_{2}.
\end{eqnarray}
That is,
\begin{eqnarray}
[d,u(n)]=-nu(n-1)\ \ \ \mbox{ for }u\in sl_{2},\; n\in \Z.
\end{eqnarray}

\bd{dyangian-module} {\em An a convention, we define a
$DY_{q}(sl_{2})$-{\em module} to be a $T(sl_{2}\otimes
\C[t,t^{-1}])$-module $W$ such that for any $w\in W$,
$$sl_{2}(n)w=0 \ \ \ \ \mbox{ for $n$ sufficiently large}$$
and such that all the defining relations of $DY_{q}(sl_{2})$
in Definition \ref{dyangian-double-algebra} hold.
A vector $w_{0}$ of a $DY_{q}(sl_{2})$-module $W$ is called
a {\em vacuum vector} if
$$u(n)w_{0}=0\ \ \ \mbox{ for }u\in sl_{2},\; n\ge 0.$$
A {\em vacuum $DY_{q}(sl_{2})$-module} is a module $W$
equipped with a vacuum vector that generates $W$.}
\ed

The following are some basic properties of a general vacuum
$DY_{q}(sl_{2})$-module:

\bl{lvacuum-property}
Let $W$ be a vacuum $DY_{q}(sl_{2})$-module with a vacuum vector $w_{0}$
as a generator. Set $F_{0}=\C w_{0}$ and $F_{k}=0$ for $k<0$.
For any positive integer $k$, we define $F_{k}$ to be the linear span of
the vectors
$$a_{1}(-m_{1})\cdots a_{r}(-m_{r})w_{0}$$
for $r\ge 1,\; a_{1},\dots, a_{r}\in sl_{2},\;
m_{1},\dots,m_{r}\ge 1$ with $m_{1}+\cdots +m_{r}\le k$.
Then the subspaces $F_{k}$ for $k\in \Z$ form an increasing filtration
of $W$ and for any $a\in sl_{2},\; m,k\in \Z$,
\begin{eqnarray}\label{eproperty}
a(m)F_{k}\subset F_{k-m}.
\end{eqnarray}
Furthermore, $T_{m}w_{0}=0$ for $m\ge 1$.
\el

\begin{proof} We first prove (\ref{eproperty}). It is true for $k<0$
as $F_{k}=0$ by definition. With $F_{0}=\C w_{0}$, we see that
(\ref{eproperty}) holds for $k=0$.
Assume $k\ge 1$. {}From definition, (\ref{eproperty}) always holds for $m<0$.
Let $a,b\in \{ e,f,h\}$. {}From the defining relations
  of $DY_{q}(sl_{2})$ we have
\begin{eqnarray}\label{erelation-dy}
a(m)b(n)=\pm b(n)a(m)+\sum_{i,j\ge 0, \;i+j\ge 1}\lambda_{ij}
b(n+i)a(m+j)+\alpha h(m+n)
\end{eqnarray}
for all $m,n\in \Z$, where $\lambda_{ij},\alpha\in \C$, depending on $a,b$.
Using this fact and induction on $k$
we obtain (\ref{eproperty}), noticing that $a(m)F_{0}=0$ for $m\ge 0$.
{}From (\ref{eproperty}) we get
$$T_{m}w_{0}= T_{m}F_{0}\subset F_{-m}=0\ \ \ \mbox{  for }m\ge 1.$$
It also follows from (\ref{eproperty}) that $\cup_{k\ge 0}F_{k}$
is a submodule of $W$.
Since $w_{0}$ generates $W$, we must have $W=\cup_{k\ge 0}F_{k}$.
This proves that the subspaces $F_{k}$ for $k\in \Z$
form an increasing filtration of $W$.
\end{proof}

\bl{cspanning-property2} Let $W$ be a vacuum $DY_{q}(sl_{2})$-module
with a vacuum vector $w_{0}$ as a generator. For $n\in \N$, define
$E_{n}$ to be the linear span of the vectors
$$u^{(1)}(m_{1})\cdots u^{(r)}(m_{r})w_{0}$$
for $0\le r\le n,\; u^{(i)}\in \{ e,f,h\},\; m_{i}\in \Z$.
Then the subspaces $E_{n}$ for $n\in \N$ form an
increasing filtration of
$W$ and for each $n\in \N$, $E_{n}$ is linearly spanned by the vectors
\begin{eqnarray}
e(-m_{1})\cdots e(-m_{r})f(-n_{1})\cdots f(-n_{s})h(-k_{1})\cdots
h(-k_{l})w_{0},
\end{eqnarray}
where $r,s,t\ge 0$ and $m_{i},n_{j},k_{t}$ are positive integers such that
$$m_{1}>\cdots >m_{r},\ \ n_{1}>\cdots >n_{s},\ \ \
k_{1}\ge \cdots \ge k_{l}, \ \ \ r+s+l\le n.$$
\el

\begin{proof} As $w_{0}$ generates $W$,
the subspaces $E_{n}$ for $n\in \N$ form an increasing filtration
for $W$. It remains to prove the spanning property. For any
nonnegative integer $n$, let $E'_{n}$ be the span of the vectors
$$a_{1}(-m_{1})\cdots a_{r}(-m_{r})w_{0}$$
for $0\le r\le n,\; a_{1},\dots, a_{r}\in sl_{2},\;
m_{1},\dots,m_{r}\ge 1$. By definition, $E_{n}'\subset E_{n}$ for
$n\ge 0$. Using induction (on $k$) and (\ref{erelation-dy}) we get
\begin{eqnarray*}
& &a(m)E'_{k}\subset E'_{k+1}\;\;\;\mbox{ if }m< 0,\\
& &a(m)E'_{k}\subset E'_{k}\;\;\;\mbox{ if }m\ge 0
\end{eqnarray*}
for any $a\in sl_{2},\; m\in \Z,\; k\in \N$. Using this and induction
we have $E_{n}\subset E_{n}'$ for $n\ge 0$.
Thus $E_{n}=E_{n}'$ for all $n\ge 0$. For
every nonnegative integer $n$, {}from Lemma \ref{lvacuum-property},
the subspaces $E_{n}\cap F_{m}$ for $m\in \N$ form an increasing
filtration  of $E_{n}$. The spanning property of $E_{n}$ follows
from this filtration and (\ref{erelation-dy}).
\end{proof}

The following is a tautological construction of a vacuum module. Let
$d$ be the derivation of $T$ such that
$$d(a\otimes t^{n})=-n(a\otimes t^{n-1})\ \ \ \mbox{ for }a\in
 sl_{2},\; n\in \Z.$$
 Set
$$T_{+}=\sum_{n\ge 1}T_{n}\;\;\;\mbox{ and }\ J=T\C[d]T_{+}.$$
With $J$ a left ideal of $T$, $T/J$ is a
(left) $T$-module and for any $v\in T$, $sl_{2}(n)(v+J)=0$ for $n$
sufficiently large.

\bd{dvq}
{\em We define $V_{q}$ to be the quotient
$T$-module of $T/J$, modulo all the defining relations of $DY_{q}(sl_{2})$.
Denote by ${\bf 1}$ the image of $1$ in $V_{q}$.}
\ed

{}From the construction, $(V_{q},{\bf 1})$
is a vacuum $DY_{q}(sl_{2})$-module.
As $dJ\subset J$, $V_{q}$ admits an action of $d$ such that
\begin{eqnarray}\label{ed-operator}
d{\bf 1}=0,\ \ \ \ [d,u(x)]=\frac{d}{dx}u(x)\ \ \ \mbox{ for }u\in sl_{2}.
\end{eqnarray}
It is clear that for any vacuum $DY_{q}(sl_{2})$-module $(W,w_{0})$
on which $d$ acts such that
$$dw_{0}=0\;\;\mbox{ and }\;\; [d,u(x)]=\frac{d}{dx}u(x)
\;\;\;\mbox{ for }u\in sl_{2},$$
there exists a unique $DY_{q}(sl_{2})$-module homomorphism
{}from $V_{q}$ to $W$, sending ${\bf 1}$ to $w_{0}$.

We are going to show that $V_{q}$ has a certain normal basis and
there is a canonical quantum vertex algebra structure on $V_{q}$. To
show that $V_{q}$ has a certain normal basis, we shall construct a
vacuum $DY_{q}(sl_{2})$-module with this property.

\br{rlie-superalgebra} {\em Here, we construct a vertex
superalgebra. Let $\g=\g_{0}\oplus \g_{1}$ be a three-dimensional
Lie superalgebra with $\g_{0}=\C \bar{h}$ (the even subspace) and
$\g_{1}=\C \bar{e}\oplus \C \bar{f}$ (the odd subspace), where
\begin{eqnarray}
[\bar{e},\bar{e}]=[\bar{f},\bar{f}]=0,\ \ \ \
[\bar{e},\bar{f}]=\bar{h},\ \ \
[\bar{h},\bar{e}]=0,\ \ \ [\bar{h},\bar{f}]=0,\ \ \ [\bar{h},\bar{h}]=0.
\end{eqnarray}
(One can show that this is indeed a Lie superalgebra by embedding
$\g$ into the Clifford algebra over the ring $\C [\bar{h}]$
associated with the vector space $\g_{1}$ equipped with a symmetric
bilinear form.) Form the loop Lie superalgebra $L(\g)=\g\otimes
\C[t,t^{-1}]$. Viewing $\C$ as a trivial $\g\otimes \C[t]$-module,
we form the induced $L(\g)$-module
$$V_{L(\g)}=U(L(\g))\otimes_{U(\g\otimes \C[t])}\C.$$
Set
$${\bf 1}=1\otimes 1\in V_{L(\g)}.$$
It follows from the P-B-W theorem that $V_{L(\g)}$
has a basis $\Bar{B}$, consisting of the vectors
\begin{eqnarray}
\bar{e}(-m_{1})\cdots \bar{e}(-m_{r})\bar{f}(-n_{1})\cdots
\bar{f}(-n_{s})\bar{h}(-k_{1})\cdots
\bar{h}(-k_{l}){\bf 1},
\end{eqnarray}
where $r,s,t\ge 0$ and $m_{i},n_{j},k_{t}$ are positive integers such that
$$m_{1}>\cdots >m_{r},\ \ n_{1}>\cdots >n_{s},\ \ \
k_{1}\ge \cdots \ge k_{l}.$$
Identify $\g$ as a subspace of $V_{L(\g)}$ through the map
$u\mapsto u(-1){\bf 1}$ for $u\in \g$.
Then there exists a (unique) vertex superalgebra structure
on $V_{L(\g)}$ with ${\bf 1}$ as the vacuum vector
and with $Y(u,x)=u(x)$ for $u=\bar{e},\bar{f},\bar{h}$.
We define a $\Z$-grading on $\g\otimes \C[t,t^{-1}]$ by
$$\wt (\g \otimes t^{n})=-n\ \ \ \mbox{ for }n\in \Z,$$
making $\g\otimes \C[t,t^{-1}]$ a $\Z$-graded Lie superalgebra.
Then $V$ is $\N$-graded with $V_{(0)}=\C {\bf 1}$ such that
$$u(m)V_{(n)}\subset V_{(n-m)}\;\;\;\mbox{ for }u\in \g,\; m,n\in \Z.$$}
\er

We are going to define a vacuum $DY_{q}(sl_{2})$-module
structure on the vertex superalgebra $V=V_{L(\g)}$.

\bl{lprepare} There exists a unique element $\Phi(t)\in \Hom (V,
V\otimes \C[t])$ such that
\begin{eqnarray}
& &\Phi(t){\bf 1}={\bf 1},\ \ \ \Phi(t)\bar{e}=\bar{e}\otimes t,\ \ \ \
\Phi(t)\bar{f}=\bar{f}\otimes t, \ \ \ \
\Phi(t)\bar{h}= \bar{h}\otimes t^{2},\\
& &\ \ \ \ \ \Phi(t)Y(v,x)=Y(\Phi(t-x)v,x)\Phi(t)\ \ \ \mbox{ for }v\in V.
\end{eqnarray}
Furthermore, we have
\begin{eqnarray}
& &\Phi(t)\bar{e}(x)=(t-x)\bar{e}(x)\Phi(t), \ \ \
\Phi(t)\bar{f}(x)=(t-x)\bar{f}(x)\Phi(t), \\
& &\Phi(t)\bar{h}(x)=(t-x)^{2}\bar{h}(x)\Phi(t),
\end{eqnarray}
and $[\D,\Phi(t)]=\frac{d}{dt}\Phi(t)$,
$\Phi(x)\Phi(t)=\Phi(t)\Phi(x)$.\el

\begin{proof} Let us equip $\C[t]$ with the  vertex algebra
structure for which  $1$ is the vacuum vector and
$$Y(p(t),x)q(t)=\left(e^{-x(d/dt)}p(t)\right) q(t)
=p(t-x)q(t)$$ for $p(t),q(t)\in \C[t]$. Then equip $V\otimes \C[t]$
with the tensor product vertex superalgebra structure where we
denote the vertex operator map by $Y_{ten}$. Thus
$$Y_{ten}(u\otimes t^{n},x)=Y(u,x)\otimes (t-x)^{n}
\ \ \ \mbox{ for }u\in V,\; n\in \Z.$$
We have
\begin{eqnarray*}
& &[Y_{ten}(\bar{e}\otimes t,x_{1}),Y_{ten}(\bar{e}\otimes t,x_{2})]_{+}
=[Y(\bar{e},x_{1}), Y(\bar{e},x_{2})]_{+}\otimes (t-x_{1})(t-x_{2})=0,\\
& &[Y_{ten}(\bar{f}\otimes t,x_{1}),Y_{ten}(\bar{f}\otimes t,x_{2})]_{+}
=[Y(\bar{f},x_{1}), Y(\bar{f},x_{2})]_{+}\otimes (t-x_{1})(t-x_{2})=0,\\
& &[Y_{ten}(\bar{h}\otimes t^{2},x_{1}),Y_{ten}(\bar{e}\otimes t,x_{2})]
=[Y(\bar{h},x_{1}), Y(\bar{e},x_{2})]\otimes (t-x_{1})^{2}(t-x_{2})=0,\\
& &[Y_{ten}(\bar{h}\otimes t^{2},x_{1}),Y_{ten}(\bar{f}\otimes t,x_{2})]
=[Y(\bar{h},x_{1}), Y(\bar{f},x_{2})]\otimes (t-x_{1})^{2}(t-x_{2})=0,\\
& &[Y_{ten}(\bar{e}\otimes t,x_{1}),Y_{ten}(\bar{f}\otimes t,x_{2})]_{+}
=[Y(\bar{e},x_{1}),Y(\bar{f},x_{2})]_{+}\otimes (t-x_{1})(t-x_{2})\\
& &\ \
=x_{1}^{-1}\delta\left(\frac{x_{2}}{x_{1}}\right)
Y(\bar{h},x_{2})\otimes (t-x_{1})(t-x_{2})
=x_{1}^{-1}\delta\left(\frac{x_{2}}{x_{1}}\right)Y(\bar{h},x_{2})\otimes
(t-x_{2})^{2},\\
& &[Y_{ten}(\bar{h}\otimes t^{2},x_{1}),Y_{ten}(\bar{h}\otimes t^{2},x_{2})]
=[Y(\bar{h},x_{1}),Y(\bar{h},x_{2})]\otimes (t-x_{1})^{2}(t-x_{2})^{2}=0.
\end{eqnarray*}
It follows that there exists a (unique) vertex-superalgebra
homomorphism $\theta$ {}from $V$ to $V\otimes \C[t]$ such that
\begin{eqnarray}
\theta (\bar{e})=\bar{e}\otimes t,\ \ \ \
\theta (\bar{f})=\bar{f}\otimes t, \ \ \ \
\theta (\bar{h})= \bar{h}\otimes t^{2}.
\end{eqnarray}
Let us alternatively denote  by $\Phi(t)$ the vertex superalgebra
homomorphism $\theta$ (from $V$ to $V\otimes \C[t]$). Then
$\Phi(t){\bf 1}={\bf 1}$, $\Phi(t)(\bar{e})=\bar{e}\otimes t$,
$\Phi(t)(\bar{f})=\bar{f}\otimes t$,
$\Phi(t)(\bar{h})=\bar{h}\otimes t^{2}$. Furthermore, for $u,v\in
V$, we have
$$\Phi(t)Y(u,x)v=\theta(Y(u,x)v)=Y_{ten}(\theta(u),x)\theta(v)
=Y(\Phi(t-x)u,x)\Phi(t)v,$$ where $Y$ is viewed as a $\C[t]$-map.
The rest follows immediately.
\end{proof}

\bp{pyangian-double-q}
Let $q$ be any nonzero complex number and let $V=V_{L(\g)}$
be the vertex superalgebra
as in Lemma \ref{lprepare}.
The assignment
\begin{eqnarray*}
e(x)=\bar{e}(x)\Phi(q+x),\ \ \
f(x)=\bar{f}(x)\Phi(-q+x),\ \ \
h(x)=q \bar{h}(x)\Phi(q+x)\Phi(-q+x)
\end{eqnarray*}
uniquely defines a vacuum $DY_{q}(sl_{2})$-module structure on $V$
with ${\bf 1}$ as the generating vacuum vector and
\begin{eqnarray}
[\D,e(x)]=\frac{d}{dx}e(x),\ \ \ [\D,f(x)]=\frac{d}{dx}f(x),\ \ \
[\D,h(x)]=\frac{d}{dx}h(x).
\end{eqnarray}
 Furthermore, for $n\in \N$, define $E_{n}$ to be the linear
span of the vectors
$$u^{(1)}(m_{1})\cdots u^{(r)}(m_{r}){\bf 1}$$
for $0\le r\le n,\; u^{(i)}\in \{ e,f,h\},\; m_{i}\in \Z$.
Then $E_{n}$ has a basis consisting of the vectors
\begin{eqnarray}
e(-m_{1})\cdots e(-m_{r})f(-n_{1})\cdots f(-n_{s})h(-k_{1})\cdots
h(-k_{l}){\bf 1},
\end{eqnarray}
where $r,s,t\ge 0$ and $m_{i},n_{j},k_{t}$ are positive integers such that
$$m_{1}>\cdots >m_{r},\ \ n_{1}>\cdots >n_{s},\ \ \
k_{1}\ge \cdots \ge k_{l},\ \ \ r+s+l\le n.$$
\ep

\begin{proof} Using Lemma \ref{lprepare} we have
\begin{eqnarray*}
e(x_{1})e(x_{2})
&=&\bar{e}(x_{1})\Phi(q+x_{1})\bar{e}(x_{2})\Phi(q+x_{2})\\
&=&\bar{e}(x_{1})\bar{e}(x_{2})\Phi(q+x_{1})\Phi(q+x_{2})(q+x_{1}-x_{2})\\
&=&-\bar{e}(x_{2})\bar{e}(x_{1})\Phi(q+x_{2})\Phi(q+x_{1})(q+x_{1}-x_{2})\\
&=&-\bar{e}(x_{2})\Phi(q+x_{2})\bar{e}(x_{1})
\Phi(q+x_{1})(q+x_{1}-x_{2})(q+x_{2}-x_{1})^{-1}\\
&=&\frac{q+x_{1}-x_{2}}{-q-x_{2}+x_{1}}
e(x_{2})e(x_{1}),\\
f(x_{1})f(x_{2})
&=&\bar{f}(x_{1})\Phi(-q+x_{1})\bar{f}(x_{2})\Phi(-q+x_{2})\\
&=&\bar{f}(x_{1})\bar{f}(x_{2})\Phi(-q+x_{1})\Phi(-q+x_{2})(-q+x_{1}-x_{2})\\
&=&-\bar{f}(x_{2})\bar{f}(x_{1})\Phi(-q+x_{2})\Phi(-q+x_{1})(-q+x_{1}-x_{2})\\
&=&-\bar{f}(x_{2})\Phi(-q+x_{2})\bar{f}(x_{1})
\Phi(-q+x_{1})(-q+x_{1}-x_{2})(-q+x_{2}-x_{1})^{-1}\\
&=&\frac{-q+x_{1}-x_{2}}{q-x_{2}+x_{1}}
f(x_{2})f(x_{1}),
\end{eqnarray*}
\begin{eqnarray*}
& &[e(x_{1}),f(x_{2})]\\
&=&\bar{e}(x_{1})\Phi(q+x_{1})\bar{f}(x_{2})\Phi(-q+x_{2})
-\bar{f}(x_{2})\Phi(-q+x_{2})\bar{e}(x_{1})\Phi(q+x_{1})\\
&=&(q+x_{1}-x_{2})\bar{e}(x_{1})\bar{f}(x_{2})\Phi(q+x_{1})\Phi(-q+x_{2})\\
& &\ \ \ \
-(-q+x_{2}-x_{1})\bar{f}(x_{2})\bar{e}(x_{1})\Phi(q+x_{1})\Phi(-q+x_{2})\\
&=&(q+x_{1}-x_{2})(\bar{e}(x_{1})\bar{f}(x_{2})+\bar{f}(x_{2})\bar{e}(x_{1}))
\Phi(q+x_{1})\Phi(-q+x_{2})\\
&=&(q+x_{1}-x_{2})x_{1}^{-1}\delta\left(\frac{x_{2}}{x_{1}}\right)
\bar{h}(x_{2})
\Phi(q+x_{1})\Phi(-q+x_{2})\\
&=&q x_{1}^{-1}\delta\left(\frac{x_{2}}{x_{1}}\right)\bar{h}(x_{2})
\Phi(q+x_{2})\Phi(-q+x_{2})\\
&=&x_{1}^{-1}\delta\left(\frac{x_{2}}{x_{1}}\right)h(x_{2}),
\end{eqnarray*}
\begin{eqnarray*}
& &[h(x_{1}),h(x_{2})]\\
&=&q^{2}\bar{h}(x_{1})\Phi(q+x_{1})\Phi(-q+x_{1})
\bar{h}(x_{2})\Phi(q+x_{2})\Phi(-q+x_{2})\\
& &\ \ -q^{2}\bar{h}(x_{2})\Phi(q+x_{2})\Phi(-q+x_{2})
\bar{h}(x_{1})\Phi(q+x_{1})\Phi(-q+x_{1})\\
&=&q^{2} (q+x_{1}-x_{2})(-q+x_{1}-x_{2})
\bar{h}(x_{1})\bar{h}(x_{2})\Phi(q+x_{1})\Phi(-q+x_{1})
\Phi(q+x_{2})\Phi(-q+x_{2})\\
& &\ \ -q^{2}(q+x_{2}-x_{1})(-q+x_{2}-x_{1})
\bar{h}(x_{2})\bar{h}(x_{1})\Phi(q+x_{1})\Phi(-q+x_{1})
\Phi(q+x_{2})\Phi(-q+x_{2})\\
&=&q^{2}(xq+_{1}-x_{2})(-q+x_{1}-x_{2}) [\bar{h}(x_{1}),\bar{h}(x_{2})]
\Phi(q+x_{1})\Phi(-q+x_{1})
\Phi(q+x_{2})\Phi(-q+x_{2})\\
&=&0,
\end{eqnarray*}
\begin{eqnarray*}
& &h(x_{1})e(x_{2})\\
&=&q\bar{h}(x_{1})\Phi(q+x_{1})\Phi(-q+x_{1})\bar{e}(x_{2})\Phi(q+x_{2})\\
&=&q \bar{h}(x_{1})\bar{e}(x_{2})
\Phi(q+x_{1})\Phi(-q+x_{1})\Phi(q+x_{2})
(q+x_{1}-x_{2})(-q+x_{1}-x_{2})\\
&=&q \bar{e}(x_{2})\bar{h}(x_{1})\Phi(q+x_{2})\Phi(q+x_{1})\Phi(-q+x_{1})
(q+x_{1}-x_{2})(-q+x_{1}-x_{2})\\
&=&\frac{q (q+x_{1}-x_{2})(-q+x_{1}-x_{2})}{(q+x_{2}-x_{1})^{2}}
\bar{e}(x_{2})\Phi(q+x_{2})\bar{h}(x_{1})\Phi(q+x_{1})\Phi(-q+x_{1})\\
&=&\frac{q+x_{1}-x_{2}}{-q-x_{2}+x_{1}}
e(x_{2})h(x_{1}),\\
& &h(x_{1})f(x_{2})\\
&=&q \bar{h}(x_{1})\Phi(q+x_{1})\Phi(-q+x_{1})\bar{f}(x_{2})\Phi(-q+x_{2})\\
&=&q \bar{h}(x_{1})\bar{f}(x_{2})
\Phi(q+x_{1})\Phi(-q+x_{1})\Phi(-q+x_{2})
(q+x_{1}-x_{2})(-q+x_{1}-x_{2})\\
&=&q \bar{f}(x_{2})\bar{h}(x_{1})\Phi(-q+x_{2})\Phi(q+x_{1})\Phi(-q+x_{1})
(q+x_{1}-x_{2})(-q+x_{1}-x_{2})\\
&=&\frac{q(q+x_{1}-x_{2})(-q+x_{1}-x_{2})}{(-q+x_{2}-x_{1})^{2}}
 \bar{f}(x_{2})\Phi(-q+x_{2})\bar{h}(x_{1})\Phi(q+x_{1})\Phi(-q+x_{1})\\
&=&\frac{-q+x_{1}-x_{2}}{q-x_{2}+x_{1}}
f(x_{2})h(x_{1}).
\end{eqnarray*}
This proves that $V$ becomes a $DY_{q}(sl_{2})$-module.
As $\Phi(x){\bf 1}={\bf 1}$, it is clear that
${\bf 1}$ is a vacuum vector for $DY_{q}(sl_{2})$.
Now it remains to prove that ${\bf 1}$
generates $V$ as an $DY_{q}(sl_{2})$-module.
Let $W$ be the $DY_{q}(sl_{2})$-submodule of $V$ generated by ${\bf 1}$.
Using Lemma \ref{lprepare} we have
$$\Phi(x_{1})e(x)=\Phi(x_{1})(\bar{e}(x)\Phi(q+x))
=\bar{e}(x)\Phi(q+x)(x_{1}-x)\Phi(x_{1})=(x_{1}-x)e(x)\Phi(x_{1}).$$
Similar relations also hold for $\bar{f}(x)$ and $\bar{h}(x)$.
As $\Phi(x){\bf 1}={\bf 1}$, by induction
we have $\Phi(x)W\subset W((x))$.
Then it follows that $W$ is stable under the actions of
$\bar{e}(n),\bar{f}(n),\bar{h}(n)$ for $n\in \Z$.
Thus $W=V$. This proves that
${\bf 1}$ generates $V$ as a $DY_{q}(sl_{2})$-module and then proves
that $V$ is a vacuum $DY_{q}(sl_{2})$-module.

Now we prove the last assertion.
With the spanning property having been established in
Lemma \ref{cspanning-property2} we only need to prove the independence.
Recall that $V=\coprod_{n\in \N}V_{(n)}$ is $\N$-graded with $V_{(0)}=\C {\bf 1}$.
For $n\in \N$, set
$$\bar{F}_{n}=V_{(0)}\oplus V_{(1)}\oplus \cdots \oplus V_{(n)}\subset V.$$
We know that $\Bar{F}_{n}$ has a basis consisting of the vectors
\begin{eqnarray*}
\bar{e}(-m_{1})\cdots \bar{e}(-m_{r})\bar{f}(-n_{1})\cdots
\bar{f}(-n_{s})\bar{h}(-k_{1})\cdots
\bar{h}(-k_{l}){\bf 1},
\end{eqnarray*}
where $r,s,t\ge 0$ and $m_{i},n_{j},k_{t}$ are positive integers such that
$$m_{1}>\cdots >m_{r},\ \ n_{1}>\cdots >n_{s},\ \ \
k_{1}\ge \cdots \ge k_{l},\ \ \ \sum m_{i}+\sum n_{j}+\sum k_{t}\le n.$$
{}From the commutation relations in Lemma \ref{lprepare},
we have
\begin{eqnarray*}
& &\Phi(t)\bar{e}(m)=(t\bar{e}(m)-\bar{e}(m+1))\Phi(t),\ \ \ \
\Phi(t)\bar{f}(m)=(t\bar{f}(m)-\bar{f}(m+1))\Phi(t),\\
& &\Phi(t)\bar{h}(m)=(t^{2}\bar{h}(m)-2t\bar{h}(m+1)+\bar{h}(m+2))\Phi(t)
\end{eqnarray*}
for $m\in \Z$. With $\Phi(t){\bf 1}={\bf 1}$, using induction we get
$$\Phi(t)w\equiv t^{m} w\ \ \ \mod \; \bar{F}_{n-1}[t]
\ \ \ \mbox{ for }w\in \bar{F}_{n},\; n\ge 0,$$
where $m$ is a nonnegative integer depending on $w$.
As
$$e(x)=\bar{e}(x)\Phi(q+x)=\sum_{j\ge 0}\frac{1}{j!}x^{j}\bar{e}(x)\Phi^{(j)}(q),$$
for any $m\in \Z$  we have
$$e(m)=\sum_{i\ge 0}\frac{1}{i!}\bar{e}(m+i)\Phi^{(i)}(q).$$
For $u\in \{ e,f,h\}$ and for $m\in \Z$, $w\in \bar{F}_{k}$,
we have
$$u(m)w\equiv \alpha \bar{u}(m)w \;\;\; \mod\; \bar{F}_{k-m-1}$$
for some nonzero complex number $\alpha$.
It follows immediately that $E_{n}$ has a basis as claimed.
\end{proof}

With $V_{q}$ being universal, from Proposition
\ref{pyangian-double-q} we immediately have:

 \bc{cnormal-basis} For $n\in \N$, let $E_{n}$
be the subspace of $V_{q}$, linearly spanned by the vectors
$$u^{(1)}(m_{1})\cdots u^{(r)}(m_{r}){\bf 1}$$
for $0\le r\le n,\; u^{(i)}\in \{ e,f,h\},\; m_{i}\in \Z$. Then the
subspaces $E_{n}$ for $n\ge 0$ form an increasing filtration of
$V_{q}$ and for each $n\ge 0$, $E_{n}$ has a basis consisting of the
vectors
\begin{eqnarray}
e(-m_{1})\cdots e(-m_{r})f(-n_{1})\cdots f(-n_{s})h(-k_{1})\cdots
h(-k_{l}){\bf 1},
\end{eqnarray}
where $r,s,t\ge 0$ and $m_{i},n_{j},k_{t}$ are positive integers
such that
$$m_{1}>\cdots >m_{r},\ \ n_{1}>\cdots >n_{s},\ \ \
k_{1}\ge \cdots \ge k_{l},\ \ \ r+s+l\le n.$$ \ec

In view of Corollary \ref{cnormal-basis}, we can and we should
consider $sl_{2}$ as a subspace of $V_{q}$ through the map $u\mapsto
u(-1){\bf 1}$ for $u\in sl_{2}$. The following is our main result:

\bt{tmain-dy-0} Let $q$ be any nonzero complex number and let
$(V_{q},{\bf 1})$ be the universal vacuum $DY_{q}(sl_{2})$-module.
There exists one and only one weak quantum vertex algebra structure
on $V_{q}$ with ${\bf 1}$ as the vacuum vector such that
$$Y(e,x)=e(x),\ \ Y(f,x)=f(x), \ \ \ Y(h,x)=h(x),$$
and the weak quantum vertex algebra $V_{q}$ is nondegenerate.
Furthermore, for any $DY_{q}(sl_{2})$-module $W$,
there exists one and only one $V_{q}$-module structure $Y_{W}$
 on $W$ such that
$$Y_{W}(e,x)=e(x),\ \ Y_{W}(f,x)=f(x),
\ \ Y_{W}(h,x)=h(x).$$
\et

\begin{proof} We shall follow the procedure outlined in
\cite{li-qva1} and \cite{li-qva2}. Let $W$ be any
$DY_{q}(sl_{2})$-module
and let $\overline{W}=V_{q}\oplus W$ be the direct sum module.
Set $U=\{ e(x),f(x),h(x)\}\subset \E(\overline{W})$. {}From the defining
relations, $U$ is an $\S$-local subset. By (\cite{li-qva1}, Theorem
5.8), $U$ generates a
weak quantum vertex algebra $V_{\overline{W}}$ where
the identity operator $1_{\overline{W}}$ is the vacuum vector
and $Y_{\E}$ denotes the vertex operator map.
Furthermore, the vector space $\overline{W}$ is a faithful
$V_{\overline{W}}$-module
with $Y_{\overline{W}}(a(x),x_{0})=a(x_{0})$ for $a(x)\in V_{\overline{W}}$.
It follows from (\cite{li-qva1}, Proposition 6.7) and the defining
relations of $DY_{q}(sl_{2})$ that $V_{\overline{W}}$ is
a vacuum $DY_{q}(sl_{2})$-module with $e(x_{0}), f(x_{0}),h(x_{0})$ acting as
$Y_{\E}(e(x),x_{0}), Y_{\E}(f(x),x_{0}), Y_{\E}(h(x),x_{0})$.
As $V_{q}$ is universal,
there exists a $DY_{q}(sl_{2})$-module homomorphism $\psi$ from $V_{q}$ to
$V_{\overline{W}}$, sending ${\bf 1}$ to $1_{\overline{W}}$.
Since $V_{q}$ is a  $DY_{q}(sl_{2})$-submodule of $\overline{W}$,
it follows that $\psi$ maps $V_{q}$ into $V_{q}\subset \overline{W}$.
Notice that $V_{q}$ as a $DY_{q}(sl_{2})$-module
is generated by ${\bf 1}$ and that we have the operator $d$
on $V_{q}$ with the property (\ref{ed-operator}). Now we can apply
Theorem 6.3 of \cite{li-qva1}, asserting that there exists one and only one
weak quantum vertex algebra structure on $V_{q}$ with the required
properties.
It follows from Theorem 6.5 of \cite{li-qva1} that
$\overline{W}$ is a $V_{q}$-module with $W$ as a submodule.

Now it remains to prove that $V_{q}$ is nondegenerate.
For $n\in \N$, define $E_{n}$ to be the linear span of the vectors
$$u^{(1)}(m_{1})\cdots u^{(r)}(m_{r}){\bf 1}$$
for $0\le r\le n,\; u^{(i)}\in \{ e,f,h\},\; m_{i}\in \Z$. By
Proposition 3.15 of \cite{li-qva2}, the subspaces $E_{n}$ $(n\in
\N)$ form an increasing filtration of $V_{q}$ with $E_{0}=\C {\bf
1}$ such that $a_{k}E_{n}\subset E_{m+n}$ for $a\in E_{m},\; m,n\in
\N,\; k\in \Z$. Denote by $Gr_{E}(V_{q})$ the associated nonlocal
vertex algebra. Notice that $e,f,h\in E_{1}$. Let
$\hat{e},\hat{f},\hat{h}$ denote the images of $e,f,h$ in
$E_{1}/E_{0}\subset Gr_{E}(V_{q})$. Then
$\{\hat{e},\hat{f},\hat{h}\}$ is a generating subset of
$Gr_{E}(V_{q})$ and we have
\begin{eqnarray*}
& &\hat{e}(x_{1})\hat{e}(x_{2})
=\frac{q+x_{1}-x_{2}}{-q+x_{1}-x_{2}}\hat{e}(x_{2})\hat{e}(x_{1}),\\
& &\hat{f}(x_{1})\hat{f}(x_{2})
=\frac{-q+x_{1}-x_{2}}{q+x_{1}-x_{2}}\hat{f}(x_{2})\hat{f}(x_{1}),\\
& &\hat{e}(x_{1})\hat{f}(x_{2})=\hat{f}(x_{2})\hat{e}(x_{1}),\\
& &\hat{h}(x_{1})\hat{h}(x_{2})=\hat{h}(x_{2})\hat{h}(x_{1}),\\
& &\hat{h}(x_{1})\hat{e}(x_{2})
=\frac{q+x_{1}-x_{2}}{-q+x_{1}-x_{2}}\hat{e}(x_{2})\hat{h}(x_{1}),\\
& &\hat{h}(x_{1})\hat{f}(x_{2})
=\frac{-q+x_{1}-x_{2}}{q+x_{1}-x_{2}}\hat{f}(x_{2})\hat{h}(x_{1}).
\end{eqnarray*}
{}From Corollary \ref{cnormal-basis}, for each $n\ge 0$,
$E_{n+1}/E_{n}$ has a basis consisting of the vectors
\begin{eqnarray*}
\hat{e}(-m_{1})\cdots \hat{e}(-m_{r})\hat{f}(-n_{1})\cdots
\hat{f}(-n_{s})\hat{h}(-k_{1})\cdots
\hat{h}(-k_{l}){\bf 1},
\end{eqnarray*}
where $r,s,t\ge 0$ and $m_{i},n_{j},k_{t}$ are positive integers such that
$$m_{1}>\cdots >m_{r},\ \ n_{1}>\cdots >n_{s},\ \ \
k_{1}\ge \cdots \ge k_{l},\ \ \ r+s+l=n+1.$$
It was proved in [KL] that $Gr_{E}(V_{q})$ is nondegenerate.
Then by (\cite{li-qva2}, Proposition 3.14), $V_{q}$ is nondegenerate.
\end{proof}

\br{rcontinue} {\em In the defining relations of the algebra
$DY_{q}(sl_{2})$, let us use the following expansion
$$\frac{1}{x_{1}-x_{2}\pm q}=\sum_{i\ge 0}(\mp q)^{i}(x_{1}-x_{2})^{-1-i}
\in \C[x_{2}][[x_{1}^{-1}]].$$ By doing this one gets a new algebra
which we denote by $DY_{q}^{\infty}(sl_{2})$. Unlike
$DY_{q}(sl_{2})$, the algebra $DY_{q}^{\infty}(sl_{2})$ only admits
modules of lowest weight type. To relate $DY_{q}^{\infty}(sl_{2})$
with quantum vertex algebras we shall need a new theory.} \er

\section{Quasi-compatibility and quasi-modules-at-infinity
for nonlocal vertex algebras}

In this section we study quasi-compatible subsets of
$\Hom(W,W((x^{-1})))$ for a general vector space $W$ and we show
that {}from any quasi-compatible subset, one can construct a
canonical nonlocal vertex algebra. We formulate a notion of quasi
module-at-infinity for a nonlocal vertex algebra and we show that
the starting vector space $W$ is naturally a quasi
module-at-infinity for the nonlocal vertex algebra generated by a
quasi-compatible subset. The theory and the results of this section
are analogous to those in \cite{li-qva1}.

Let $W$ be any vector space over $\C$, which is fixed throughout this section.
Set
\begin{eqnarray}
\E^{o}(W)=\Hom (W,W((x^{-1})))\subset (\End W)[[x,x^{-1}]].
\end{eqnarray}
Denote by $1_{W}$ the identity operator on $W$, a distinguished
element of $\E^{o}(W)$. Note that $\E^{o}(W)$ is naturally a vector
space over the field $\C((x^{-1}))$. Let $G$ denote the group of
linear transformations on $\C$:
\begin{eqnarray}
G=\{g(z)=c_{0}z+c_{1}\;|\; c_{0}\in \C^{\times},\; c_{1}\in \C\}.
\end{eqnarray}
Group $G$ acts on $\E^{o}(W)$ with $g\in G$ acting as $R_{g}$
defined by
\begin{eqnarray}
R_{g}a(x)=a(g(x))\ \ \ \mbox{ for }a(x)\in \E^{o}(W),
\end{eqnarray}
where as a convention
$$a(g(x))=\sum_{n\in \Z}a_{n}(c_{0}x+c_{1})^{-n-1}
=\sum_{n\in \Z}\sum_{i\in \N}\binom{-n-1}{i}
c_{0}^{-n-1-i}c_{1}^{i}a_{n}x^{-n-1-i}$$
for $a(x)=\sum_{n\in \Z}a_{n}x^{-n-1}$.

\bd{dcompatibility}
{\em An ordered sequence $a_{1}(x),\dots, a_{r}(x)$ in $\E^{o}(W)$
is said to be
{\em quasi-compatible} if there exists a nonzero polynomial
$p(x_{1},x_{2})$ such that
\begin{eqnarray}\label{equasicompatible-condition-sequence}
\left(\prod_{1\le i<j\le r}p(x_{i},x_{j})\right)
a_{1}(x_{1})a_{2}(x_{2})\cdots a_{r}(x_{r})
\in \Hom (W,W((x_{1}^{-1},x_{2}^{-1},\dots,x_{r}^{-1}))).
\end{eqnarray}
A subset $S$ of $\E^{o}(W)$ is said to be {\em quasi-compatible} if
every finite sequence in $S$ is quasi-compatible.
A sequence $a_{1}(x),\dots, a_{r}(x)$ is said to be
{\em compatible} if there exists a nonnegative integer $k$ such that
\begin{eqnarray}\label{ecompatible-condition-sequence}
\left(\prod_{1\le i<j\le r}(x_{i}-x_{j})^{k}\right)
a_{1}(x_{1})a_{2}(x_{2})\cdots a_{r}(x_{r})
\in \Hom (W,W((x_{1}^{-1},x_{2}^{-1},\dots,x_{r}^{-1}))).
\end{eqnarray}
A subset $S$ of $\E^{o}(W)$ is said to be {\em compatible} if
every finite sequence in $S$ is compatible.}
\ed

\br{rexistence} {\em  Note that for any $A(x_{1},x_{2})\in \Hom
(W,W((x_{1}^{-1},x_{2}^{-1})))$, $A(x_{2},x_{2})$ exists in $\Hom
(W,W((x_{2}^{-1})))$ and $A(x_{2},x_{2})=0$ if and only if
  $A(x_{1},x_{2})=0$. Furthermore,
$$A(x_{2}+x_{0},x_{2})
=\left(e^{x_{0}\frac{\partial}{\partial x_{1}}}
A(x_{1},x_{2})\right)|_{x_{1}=x_{2}}$$
exists in $(\Hom (W,W((x_{2}^{-1}))))[[x_{0}]]$ and
$A(x_{2}+x_{0},x_{2})\ne 0$ if $A(x_{1},x_{2})\ne 0$.
We also have
\begin{eqnarray*}
& &x_{2}^{-1}\delta\left(\frac{x_{1}}{x_{2}}\right)A(x_{1},x_{2})
=x_{2}^{-1}\delta\left(\frac{x_{1}}{x_{2}}\right)A(x_{2},x_{2}),\\
& &\Res_{x_{1}}x_{2}^{-1}\delta\left(\frac{x_{1}}{x_{2}}\right)A(x_{1},x_{2})
=A(x_{2},x_{2}),\\
& &\Res_{x_{1}}x_{1}^{-1}\delta\left(\frac{x_{2}+x_{0}}{x_{1}}\right)
A(x_{1},x_{2})
=A(x_{2}+x_{0},x_{2}).
\end{eqnarray*}}
\er

Let $\C(x_{1},x_{2})$ denote the field of rational functions.
We have fields $\C((x_{1}^{-1}))((x_{2}^{-1}))$ and
$\C((x_{1}^{-1}))((x_{2}))$ of formal series.
As these fields
contain $\C[x_{1},x_{2}]$ as a subring,
there exist unique field-embeddings
\begin{eqnarray}
\iota_{x_{1},\infty;x_{2},\infty}:& & \C(x_{1},x_{2})\rightarrow
\C((x_{1}^{-1}))((x_{2}^{-1})),\\
\iota_{x_{1},\infty;x_{2},0}:& & \C(x_{1},x_{2})\rightarrow
\C((x_{1}^{-1}))((x_{2})).
\end{eqnarray}
Let $\C(x)$ denote the field of rational functions. Define
field-embeddings
$$\iota_{x,0}: \C(x)\rightarrow \C((x)),$$
 sending $f(x)$ to the
formal Laurent series expansion of $f(x)$ at $x=0$, and
$$\iota_{x,\infty}: \C(x)\rightarrow \C((x^{-1})),$$
 sending $f(x)$ to
the formal Laurent series expansion of $f(x)$ at $x=\infty$.

\br{rcancellation} {\em We shall often use the following simple
fact. Suppose that $U$ is a vector space over $\C$ and $F,G\in
U((x^{-1}))((x_{0}))$ satisfy the relation
$$q(x,x_{0})F(x,x_{0})=q(x,x_{0})G(x,x_{0})$$ for some nonzero
polynomial $q(x,x_{0})$. Then $F=G$. This is simply because
$U((x^{-1}))((x_{0}))$ is a vector space over the field
$\C((x^{-1}))((x_{0}))$ and $0\ne q(x,x_{0})\in \C[x,x_{0}]\subset
\C((x^{-1}))((x_{0}))$. } \er

Note that for a quasi-compatible pair
$(a(x),b(x))$ in $\E^{o}(W)$, by definition
there exists a nonzero polynomial $p(x_{1},x_{2})$ such that
\begin{eqnarray}\label{ecompatible-condition}
p(x_{1},x_{2})a(x_{1})b(x_{2})\in \Hom (W,W((x_{1}^{-1},x_{2}^{-1}))).
\end{eqnarray}
As
\begin{eqnarray*}
& &\iota_{x,\infty;x_{0},0}\left(1/p(x_{0}+x,x)\right)\in
  \C((x^{-1}))((x_{0})),\\
& &\left( p(x_{1},x)a(x_{1})b(x)\right)|_{x_{1}=x+x_{0}}
\in (\Hom(W,W((x^{-1}))))[[x_{0}]],
\end{eqnarray*}
we have
$$\iota_{x,\infty;x_{0},0}\left(1/p(x_{0}+x,x)\right)
\left( p(x_{1},x)a(x_{1})b(x)\right)|_{x_{1}=x+x_{0}}
\in (\Hom (W,W((x^{-1}))))((x_{0})).$$

\bd{doperation-same-order}
{\em Let $(a(x),b(x))$ be a quasi-compatible pair in $\E^{o}(W)$.
We define
\begin{eqnarray*}
Y_{\E^{o}}(a(x),x_{0})b(x)=\iota_{x,\infty;x_{0},0}\left(1/p(x_{0}+x,x)\right)
\left( p(x_{1},x)a(x_{1})b(x)\right)|_{x_{1}=x+x_{0}}\in \E^{o}(W)((x_{0})),
\end{eqnarray*}
where $p(x_{1},x_{2})$ is any nonzero polynomial such that
(\ref{ecompatible-condition}) holds.}
\ed

It is easy to show that $Y_{\E^{o}}(a(x),x_{0})b(x)$ is well defined, i.e.,
the expression on the right hand side
does not depend on the choice of polynomial $p(x_{1},x_{2})$.

Write
\begin{eqnarray}
Y_{\E^{o}}(a(x),x_{0})b(x)=\sum_{n\in \Z}a(x)_{n}b(x) x_{0}^{-n-1}.
\end{eqnarray}

The following is an immediate consequence:

\bl{ltruncation} Let $(a(x),b(x))$ be a quasi-compatible pair in
$\E^{o}(W)$. Then $a(x)_{n}b(x)\in \E^{o}(W)$ for $n\in \Z$.
Furthermore, let $p(x_{1},x_{2})$ be a nonzero polynomial such that
$$p(x_{1},x_{2})a(x_{1})b(x_{2})\in \Hom (W,W((x_{1}^{-1},x_{2}^{-1})))$$
and let $k$ be an integer such that
$$x_{0}^{k}\iota_{x,\infty;x_{0},0}\left(1/p(x_{0}+x,x)\right)
\in \C((x^{-1}))[[x_{0}]].$$
Then
\begin{eqnarray}
a(x)_{n}b(x)=0\;\;\;\mbox{ for }n\ge k.
\end{eqnarray}
\el

We shall need the following result:

\bl{lproof-need}
Let $(a_{i}(x),b_{i}(x))$ $(i=1,\dots,n)$ be quasi-compatible
ordered pairs in $\E^{o}(W)$. Suppose that
\begin{eqnarray}\label{esum-lemma-case1}
\sum_{i=1}^{n}g_{i}(z,x)a_{i}(z)b_{i}(x)\in \Hom (W,W((z^{-1},x^{-1})))
\end{eqnarray}
for some polynomials $g_{1}(z,x),\dots,g_{n}(z,x)$. Then
\begin{eqnarray}\label{eassoc-sum}
\sum_{i=1}^{n}g_{i}(x+x_{0},x)Y_{\E^{o}}(a_{i}(x),x_{0})b_{i}(x)
=\left(\sum_{i=1}^{n}g_{i}(z,x)a_{i}(z)b_{i}(x)\right)|_{z=x+x_{0}}.
\end{eqnarray}
\el

\begin{proof} Let $g(z,x)$ be a nonzero polynomial such that
$$g(z,x)a_{i}(z)b_{i}(x)\in \Hom (W,W((z^{-1},x^{-1})))
\;\;\;\mbox{ for }i=1,\dots,n.$$
{}From Definition \ref{doperation-same-order}, we have
$$g(x+x_{0},x)Y_{\E^{o}}(a_{i}(x),x_{0})b_{i}(x)
=\left(g(z,x)a_{i}(z)b_{i}(x)\right)\mid_{z=x+x_{0}}
\;\;\;\mbox{ for }i=1,\dots,n.$$
Then using (\ref{esum-lemma-case1}) we have
\begin{eqnarray}
& &g(x+x_{0},x)\sum_{i=1}^{n}g_{i}(x+x_{0},x)Y_{\E^{o}}(a_{i}(x),x_{0})b_{i}(x)
\nonumber\\
&=&\sum_{i=1}^{n}g_{i}(x+x_{0},x)\left(g(z,x)a_{i}(z)b_{i}(x)\right)|_{z=x+x_{0}}
\nonumber\\
&=&\left(g(z,x)\sum_{i=1}^{n}g_{i}(z,x)a_{i}(z)b_{i}(x)\right)|_{z=x+x_{0}}
\nonumber\\
&=&g(x+x_{0},x)\left(\sum_{i=1}^{n}g_{i}(z,x)a_{i}(z)b_{i}(x)\right)|_{z=x+x_{0}}.
\end{eqnarray}
As both $\sum_{i=1}^{n}g_{i}(x+x_{0},x)
Y_{\E^{o}}(a_{i}(x),x_{0})b_{i}(x)$ and
$$\left(\sum_{i=1}^{n}g_{i}(z,x)a_{i}(z)b_{i}(x)\right)|_{z=x+x_{0}}$$
lie in $(\Hom
(W,W((x^{-1}))))((x_{0}))$, {}from Remark \ref{rcancellation},
(\ref{eassoc-sum}) follows immediately.
\end{proof}

A quasi-compatible subspace $U$ of $\E^{o}(W)$ is said to be {\em closed} if
\begin{eqnarray}
a(x)_{n}b(x)\in U\;\;\;\mbox{ for }a(x),b(x)\in U,\; n\in \Z.
\end{eqnarray}
We are going to prove that any closed quasi-compatible subspace
containing $1_{W}$ of $\E^{o}(W)$
is a nonlocal vertex algebra.
First we prove the following result:

\bl{lclosed}
Let $V$ be a closed quasi-compatible subspace of $\E^{o}(W)$.
Let $\psi(x),\phi(x),\theta(x)\in V$ and
let $f(x,y)$ be a nonzero polynomial such that
\begin{eqnarray}
& &f(x,y)\phi(x)\theta(y)\in \Hom (W,W((x^{-1},y^{-1}))),
\label{e4.31}\\
& &f(x,y)f(x,z)f(y,z)
\psi(x)\phi(y)\theta(z)\in \Hom (W,W((x^{-1},y^{-1},z^{-1}))).\label{e4.32}
\end{eqnarray}
Then
\begin{eqnarray}
& &f(x+x_{1},x)f(x+x_{2},x)f(x+x_{1},x+x_{2})
Y_{\E^{o}}(\psi(x),x_{1})Y_{\E^{o}}(\phi(x),x_{2})\theta(x)\nonumber\\
&=&\left(f(y,x)f(z,x)f(y,z)
\psi(y)\phi (z)\theta(x)\right)|_{y=x+x_{1},z=x+x_{2}}.
\end{eqnarray}
\el

\begin{proof} With (\ref{e4.31}), {}from
Definition \ref{doperation-same-order} we have
\begin{eqnarray}\label{ehphi-theta}
f(x+x_{2},x)Y_{\E^{o}}(\phi (x),x_{2})\theta (x)
=\left(f(z,x)\phi(z)\theta(x)\right)|_{z=x+x_{2}},
\end{eqnarray}
which gives
\begin{eqnarray}
& &f(y,x)f(y,x+x_{2})f(x+x_{2},x)\psi(y)
Y_{\E^{o}}(\phi (x),x_{2})\theta(x)\nonumber\\
&=&\left(f(y,x)f(y,z) f(z,x)
\psi(y)\phi(z)\theta(x)\right)|_{z=x+x_{2}}.
\end{eqnarray}
{}From (\ref{e4.32})
the expression on the right-hand side lies in
$(\Hom (W,W((y^{-1}, x^{-1}))))[[x_{2}]]$,
so does the expression on the left-hand side.
That is,
\begin{eqnarray*}
f(y,x)f(y,x+x_{2})f(x+x_{2},x)\psi(y)
Y_{\E^{o}}(\phi(x),x_{2})\theta(x)
\in (\Hom (W,W((y^{-1},x^{-1}))))[[x_{2}]].
\end{eqnarray*}
Multiplying by $\iota_{x,\infty;x_{2},0}(f(x+x_{2},x)^{-1})$, which
lies in $\C ((x^{-1}))((x_{2}))$, we have
\begin{eqnarray}
f(y,x)f(y,x+x_{2})
\psi(y)Y_{\E^{o}}(\phi(x),x_{2})\theta(x)
\in (\Hom (W,W((y^{-1},x^{-1}))))((x_{2})).
\end{eqnarray}
In view of Lemma \ref{lproof-need},
by considering the coefficient of each power of $x_{2}$, we have
\begin{eqnarray}
& &f(x+x_{1},x)f(x+x_{1},x+x_{2})
Y_{\E^{o}}(\psi(x),x_{1})Y_{\E^{o}}(\phi(x),x_{2})\theta(x)
\nonumber\\
&=&\left(f(y,x)f(y,x+x_{2})
\psi(y)(Y_{\E^{o}}(\phi(x),x_{2})\theta(x)\right)|_{y=x+x_{1}}.
\end{eqnarray}
Using this and (\ref{ehphi-theta}) we have
\begin{eqnarray*}
& &f(x+x_{1},x)f(x+x_{2},x)f(x+x_{1},x+x_{2})
Y_{\E^{o}}(\psi(x),x_{1})Y_{\E^{o}}(\phi(x),x_{2})\theta(x)\nonumber\\
&=&\left(f(y,x)f(x+x_{2},x)f(y,x+x_{2})
\psi(y)Y_{\E^{o}}(\phi (x),x_{2})\theta(x)\right)|_{y=x+x_{1}}\nonumber\\
&=&\left(f(y,x)f(z,x)f(y,z)
\psi(y)\phi (z)\theta(x)\right)|_{y=x+x_{1},z=x+x_{2}},
\end{eqnarray*}
as desired.
\end{proof}

To state our first result we shall need a new notion.

\bd{dmodule-infinity}
{\em Let $V$ be a nonlocal vertex algebra.
A (left) {\em quasi $V$-module-at-infinity} is a vector space $W$
equipped with a linear map
$$Y_{W}: V\rightarrow \Hom (W,W((x^{-1})))\subset (\End W)[[x,x^{-1}]],$$
satisfying the condition that $Y_{W}({\bf 1},x)=1_{W}$ and that for
any $u,v\in V$, there exists a nonzero polynomial $p(x_{1},x_{2})$
such that
\begin{eqnarray}
p(x_{1},x_{2})Y_{W}(u,x_{1})Y_{W}(v,x_{2})
\in \Hom (W,W((x_{1}^{-1},x_{2}^{-1})))
\end{eqnarray}
and
\begin{eqnarray}
p(x_{0}+x_{2},x_{2})Y_{W}(Y(u,x_{0})v,x_{2})=\left(p(x_{1},x_{2})
Y_{W}(u,x_{1})Y_{W}(v,x_{2})\right)|_{x_{1}=x_{2}+x_{0}}.
\end{eqnarray}
A quasi $V$-module at infinity $(W,Y_{W})$ is called a (left)
{\em $V$-module-at-infinity}
if for any $u,v\in V$, there exists a nonnegative integer $k$ such that
\begin{eqnarray}
(x_{1}-x_{2})^{k}Y_{W}(u,x_{1})Y_{W}(v,x_{2})
\in \Hom (W,W((x_{1}^{-1},x_{2}^{-1})))
\end{eqnarray}
and
\begin{eqnarray}
x_{0}^{k}Y_{W}(Y(u,x_{0})v,x_{2})=\left((x_{1}-x_{2})^{k}
Y_{W}(u,x_{1})Y_{W}(v,x_{2})\right)|_{x_{1}=x_{2}+x_{0}}.
\end{eqnarray}}
\ed

Here we make a new notation for convenience.
Let $a(x)=\sum_{n\in \Z}a_{n}x^{-n-1}$ be any formal series
(with coefficients $a_{n}$ in some vector space).
For any $m\in \Z$, we set
\begin{eqnarray}
a(x)_{\ge m}=\sum_{n\ge m}a_{n}x^{-n-1}.
\end{eqnarray}
Then for any polynomial $q(x)$ we have
\begin{eqnarray}\label{esimplefactresidule}
\Res_{x}x^{m}q(x)a(x)=\Res_{x}x^{m}q(x)a(x)_{\ge m}.
\end{eqnarray}
Now we are in a position to prove our first key result:

\bt{tclosed}
Let $V$ be a closed quasi-compatible subspace of $\E^{o}(W)$,
containing $1_{W}$.
Then $(V,Y_{\E^{o}},1_{W})$ carries the structure of a
nonlocal vertex algebra with $W$ as a faithful (left)
quasi module-at-infinity
where the vertex operator map $Y_{W}$ is given by
$Y_{W}(\alpha(x),x_{0})=\alpha(x_{0})$. Furthermore, if $V$ is
compatible, $W$ is a $V$-module-at-infinity.
\et

\begin{proof} For any $a(x)\in \E^{o}(W)$, as
$1_{W}a(x_{2})=a(x_{2})\in \Hom (W,W((x_{1}^{-1},x_{2}^{-1})))$,
by definition we have
\begin{eqnarray*}
& &Y_{\E^{o}}(1_{W},x_{0})a(x)=1_{W}a(x)=a(x),\\
& &Y_{\E^{o}}(a(x),x_{0})1_{W}=(a(x_{1})1_{W})|_{x_{1}=x+x_{0}}=a(x+x_{0})
=e^{x_{0}\frac{d}{dx}}a(x).
\end{eqnarray*}
For the assertion on the nonlocal vertex algebra structure,
it remains to prove the weak associativity, i.e.,
for $\psi,\phi,\theta\in V$, there exists a nonnegative
integer $k$ such that
\begin{eqnarray}\label{eweakassocmainthem}
(x_{0}+x_{2})^{k}Y_{\E^{o}}(\psi,x_{0}+x_{2})Y_{\E^{o}}(\phi,x_{2})\theta
=(x_{0}+x_{2})^{k}Y_{\E^{o}}(Y_{\E^{o}}(\psi,x_{0})\phi,x_{2})\theta.
\end{eqnarray}
Let $f(x,y)$ be a nonzero polynomial such that
\begin{eqnarray*}
& &f(x,y)\psi(x)\phi(y)\in \Hom (W,W((x^{-1},y^{-1}))),\\
& &f(x,y)\phi(x)\theta(y)\in \Hom (W,W((x^{-1},y^{-1}))),\\
& &f(x,y)f(x,z)f(y,z)
\psi(x)\phi(y)\theta(z)\in \Hom (W,W((x^{-1},y^{-1},z^{-1}))).
\end{eqnarray*}
By Lemma \ref{lclosed}, we have
\begin{eqnarray}\label{e5.76}
& &f(x+x_{2},x)f(x+x_{0}+x_{2},x)f(x+x_{0}+x_{2},x+x_{2})
Y_{\E^{o}}(\psi(x),x_{0}+x_{2})Y_{\E^{o}}(\phi(x),x_{2})\theta(x)
\nonumber\\
& &=\left(f(z,x)f(y,x)f(y,z)
\psi(y)\phi(z)\theta(x)\right)|_{y=x+x_{0}+x_{2},z=x+x_{2}}.
\end{eqnarray}

On the other hand, let $n\in \Z$ be {\em arbitrarily fixed}.
Since $\psi(x)_{m}\phi(x)=0$ for $m$ sufficiently large, there exists
a nonzero polynomial $p(x,y)$, depending on $n$, such that
\begin{eqnarray}\label{eges}
p(x+x_{2},x)
(Y_{\E^{o}}(\psi(x)_{m}\phi(x), x_{2})\theta(x)
=\left(p(z,x)(\psi(z)_{m}\phi(z))\theta(x)\right)|_{z=x+x_{2}}
\end{eqnarray}
for {\em all} $m\ge n$.
With $f(x,y)\psi(x)\phi(y)\in \Hom (W,W((x^{-1},y^{-1})))$,
{}from Definition \ref{doperation-same-order} we have
\begin{eqnarray}\label{ethis}
f(x_{2}+x_{0},x_{2})(Y_{\E^{o}}(\psi(x_{2}),x_{0})\phi(x_{2}))\theta(x)
=\left(f(y,x_{2})\psi(y)\phi(x_{2})\theta(x)\right)|_{y=x_{2}+x_{0}}.
\end{eqnarray}
Using (\ref{esimplefactresidule}), (\ref{eges}) and (\ref{ethis}) we get
\begin{eqnarray}\label{e5.78}
& &\Res_{x_{0}}x_{0}^{n}f(x+x_{0}+x_{2},x)f(x+x_{0}+x_{2},x+x_{2})
p(x+x_{2},x)\nonumber\\
& &\ \ \ \ \cdot Y_{\E^{o}}(Y_{\E^{o}}(\psi(x),x_{0})\phi(x), x_{2})\theta(x)\nonumber\\
&=&\Res_{x_{0}}x_{0}^{n}f(x+x_{0}+x_{2},x)
f(x+x_{0}+x_{2},x+x_{2})p(x+x_{2},x)\nonumber\\
& &\ \ \ \ \cdot
Y_{\cal{E}}(Y_{\E^{o}}(\psi(x),x_{0})_{\ge n}\phi(x), x_{2})\theta(x)\nonumber\\
&=&\Res_{x_{0}}x_{0}^{n}
f(x+x_{0}+x_{2},x)f(x+x_{0}+x_{2},x+x_{2})\nonumber\\
& &\ \ \ \ \cdot \left(p(z,x)
Y_{\E^{o}}(\psi(z),x_{0})_{\ge n}\phi(z))\theta(x)\right)|_{z=x+x_{2}}\nonumber\\
&=&\Res_{x_{0}}x_{0}^{n}
f(x+x_{0}+x_{2},x)f(x+x_{0}+x_{2},x+x_{2})
\nonumber\\
& &\ \ \ \ \cdot \left(p(z,x)Y_{\E^{o}}(\psi(z),x_{0})\phi(z))\theta(x)\right)|_{z=x+x_{2}}\nonumber\\
&=&\Res_{x_{0}}x_{0}^{n}\left(f(z+x_{0},x)f(z+x_{0},z)
p(z,x)(Y_{\E^{o}}(\psi(z),x_{0})\phi(z))\theta(x)\right)|_{z=x+x_{2}}\nonumber\\
&=&\Res_{x_{0}}x_{0}^{n}\left(f(y,x)f(y,z)
p(z,x)\psi(y)\phi(z)\theta(x)\right)|_{y=z+x_{0},z=x+x_{2}}.
\end{eqnarray}
Combining (\ref{e5.78}) with (\ref{e5.76}) we get
\begin{eqnarray}\label{e5.79}
& &\Res_{x_{0}}x_{0}^{n}
f(x+x_{2},x)f(x+x_{0}+x_{2},x)f(x+x_{0}+x_{2},x+x_{2})
\nonumber\\
& &\ \ \ \ \cdot
p(x+x_{2},x)Y_{\E^{o}}(\psi(x),x_{0}+x_{2})
Y_{\E^{o}}(\phi(x),x_{2})\theta (x)\nonumber\\
&=&\Res_{x_{0}}x_{0}^{n}f(x_{2}+x,x)f(x+x_{0}+x_{2},x)f(x+x_{0}+x_{2},x+x_{2})
\nonumber\\
& &\ \ \ \ \cdot p(x+x_{2},x)
Y_{\E^{o}}(Y_{\E^{o}}(\psi(x),x_{0})\phi(x), x_{2})\theta(x).
\end{eqnarray}
Notice that both sides of (\ref{e5.79}) involve only finitely
many negative powers of $x_{2}$. In view of Remark \ref{rcancellation}
we can multiply both sides by
$\iota_{x,\infty;x_{2},0}(p(x+x_{2},x)^{-1}f(x+x_{2},x)^{-1})$
(in $\C((x^{-1}))((x_{2}))$ to get
\begin{eqnarray*}
& &\Res_{x_{0}}x_{0}^{n}f(x+x_{0}+x_{2},x)
f(x+x_{0}+x_{2},x+x_{2})
Y_{\E^{o}}(\psi(x),x_{0}+x_{2})Y_{\E^{o}}(\phi(x),x_{2})\theta(x)\nonumber\\
& &=\Res_{x_{0}}x_{0}^{n}
f(x+x_{0}+x_{2},x)f(x+x_{0}+x_{2},x+x_{2})
Y_{\E^{o}}(Y_{\E^{o}}(\psi(x),x_{0})\phi(x), x_{2})\theta(x).\ \ \ \
\end{eqnarray*}
Since {\em $f(x,y)$ does not depend on $n$
and since $n$ is arbitrary}, we have
\begin{eqnarray}\label{enearfinal-new}
& &f(x+x_{0}+x_{2},x)f(x+x_{0}+x_{2},x+x_{2})
Y_{\E^{o}}(\psi(x),x_{0}+x_{2})Y_{\E^{o}}(\phi(x),x_{2})\theta(x)\nonumber\\
&=&f(x+x_{0}+x_{2},x)f(x+x_{0}+x_{2},x+x_{2})
Y_{\E^{o}}(Y_{\E^{o}}(\psi(x),x_{0})\phi(x), x_{2})\theta(x).
\end{eqnarray}
Write $f(x,y)=(x-y)^{k}g(x,y)$ for some $k\in \N,\; g(x,y)\in \C[x,y]$
with $g(x,x)\ne 0$.
Then
\begin{eqnarray*}
& &f(x+x_{0}+x_{2},x)=(x_{0}+x_{2})^{k}g(x+x_{0}+x_{2},x),\\
& &f(x+x_{0}+x_{2},x+x_{2})=x_{0}^{k}g(x+x_{0}+x_{2},x+x_{2}).
\end{eqnarray*}
Since $g(x,x)\ne 0$, we have
$$\iota_{x,\infty;z,0}g(x+z,x)^{-1}\in \C((x^{-1}))[[z]],$$
so that
$$\iota_{x,\infty;z,0}g(x+z,x)^{-1}|_{z=x_{0}+x_{2}},
\ \ \ \iota_{z,\infty;x_{0},0}g(z+x_{0},z)^{-1}|_{z=x+x_{2}}
\in \C((x^{-1}))[[x_{0},x_{2}]].$$
By cancellation, from (\ref{enearfinal-new}) we obtain
\begin{eqnarray*}
& &(x_{0}+x_{2})^{k}Y_{\E^{o}}(\psi(x),x_{0}+x_{2})Y_{\E^{o}}(\phi(x),x_{2})
\theta(x)\\
&=&(x_{0}+x_{2})^{k}
Y_{\E^{o}}(Y_{\E^{o}}(\psi(x),x_{0})\phi(x),x_{2})\theta(x),
\end{eqnarray*}
as desired. This proves that $(V,Y_{\E^{o}},1_{W})$ carries the structure
of a nonlocal vertex algebra.

Next, we prove that $W$ is a quasi module-at-infinity.
For $a(x),b(x)\in V$, there exists a nonzero polynomial $h(x,y)$
such that
$$h(x,y)a(x)b(y)\in \Hom (W,W((x^{-1},y^{-1}))).$$
Then
$$h(x_{1},x_{2})Y_{W}(a(x),x_{1})Y_{W}(b(x),x_{2})
=h(x_{1},x_{2})a(x_{1})b(x_{2})
\in \Hom (W,W((x_{1}^{-1},x_{2}^{-1})))$$
and
\begin{eqnarray*}
& &h(x_{0}+x_{2},x_{2})Y_{W}(Y_{\E^{o}}(a(x),x_{0})b(x),x_{2})\nonumber\\
&=&h(x_{0}+x_{2},x_{2})(Y_{\E^{o}}(a(x),x_{0})b(x))|_{x=x_{2}}\nonumber\\
&=&\left(h(x_{1},x_{2})
a(x_{1})b(x_{2})\right)|_{x_{1}=x_{2}+x_{0}}\nonumber\\
&=&\left(h(x_{1},x_{2})
Y_{W}(a(x),x_{1})Y_{W}(b(x),x_{2})\right)|_{x_{1}=x_{2}+x_{0}}.
\end{eqnarray*}
Therefore $W$ is a (left) quasi $V$-module-at-infinity with
$Y_{W}(\alpha(x),x_{0})=\alpha(x_{0})$ for $\alpha(x)\in V$.
Finally, if $V$ is compatible, the polynomial $h(x,y)$ is of the
form $(x-y)^{k}$ with $k\in \N$. Then $W$ is a
$V$-module-at-infinity, instead of a quasi $V$-module-at-infinity.
\end{proof}

In practice, we are often given an unnecessarily closed quasi-compatible
subspace. Next, we are going to show that every quasi-compatible
subset is contained in some closed quasi-compatible subspace.

The following is an analogue of a result in \cite{li-g1} and
\cite{li-qva1}:

\bp{pgeneratingcomplicatedone} Let $\psi_{1}(x),\dots,\psi_{r}(x),
a(x),b(x),\phi_{1}(x), \dots,\phi_{s}(x) \in \E^{o}(W)$. Assume that
the ordered sequences $(a(x), b(x))$ and
$(\psi_{1}(x),\dots,\psi_{r}(x), a(x),b(x),\phi_{1}(x),
\dots,\phi_{s}(x))$
 are quasi-compatible (compatible).
Then for any $n\in \Z$, the ordered sequence
$$(\psi_{1}(x),\dots,\psi_{r}(x),a(x)_{n}b(x),\phi_{1}(x),\dots,\phi_{s}(x))$$
 is quasi-compatible (compatible).
\ep

\begin{proof}
Let $f(x,y)$ be a nonzero polynomial such that
$$f(x,y)a(x)b(y)\in \Hom (W,W((x^{-1},y^{-1})))$$
and
\begin{eqnarray}\label{elong-exp}
& &\left(\prod_{1\le i<j\le r}f(y_{i},y_{j})\right)
\left(\prod_{1\le i\le r, 1\le j\le s}f(y_{i},z_{j})\right)
\left(\prod_{1\le i<j\le s}f(z_{i},z_{j})\right)\nonumber\\
& &\;\;\cdot f(x_{1},x_{2})
\left(\prod_{i=1}^{r}f(x_{1},y_{i})f(x_{2},y_{i})\right)
\left(\prod_{i=1}^{s}f(x_{1},z_{i})f(x_{2},z_{i})\right)\nonumber\\
& &\;\;\cdot \psi_{1}(y_{1})\cdots \psi_{r}(y_{r})
a(x_{1})b(x_{2})\phi_{1}(z_{1})\cdots \phi_{s}(z_{s})\nonumber\\
& &\in \Hom (W,W((y_{1}^{-1},\dots, y_{r}^{-1},x_{1}^{-1},x_{2}^{-1},z_{1}^{-1},\dots,z_{s}^{-1}))).
\end{eqnarray}
Set
$$P=\prod_{1\le i<j\le r}f(y_{i},y_{j}),\;\;\;\;
Q=\prod_{1\le i<j\le s}f(z_{i},z_{j}),\;\;\;\;
R=\prod_{1\le i\le r,\; 1\le j\le s}f(y_{i},z_{j}).$$
Let $n\in \Z$ be {\em arbitrarily fixed}.
There exists a nonnegative integer $k$ such that
\begin{eqnarray}\label{etruncationpsiphi}
x_{0}^{k+n}\iota_{x_{2},\infty;x_{0},0}\left(f(x_{0}+x_{2},x_{2})^{-1}\right)
\in \C ((x_{2}^{-1}))[[x_{0}]].
\end{eqnarray}
Using (\ref{etruncationpsiphi}) and Definition \ref{doperation-same-order}
we obtain
\begin{eqnarray}\label{ecompatibilitythreeproof}
& &\prod_{i=1}^{r}f(x_{2},y_{i})^{k}
\prod_{j=1}^{s}f(x_{2},z_{j})^{k}\nonumber\\
& &\ \ \ \ \cdot \psi_{1}(y_{1})\cdots\psi_{r}(y_{r})
(a(x_{2})_{n}b(x_{2}))\phi_{1}(z_{1})\cdots\phi_{s}(z_{s})\nonumber\\
&=&\Res_{x_{0}}x_{0}^{n}\prod_{i=1}^{r}f(x_{2},y_{i})^{k}
\prod_{j=1}^{s}f(x_{2},z_{j})^{k}\nonumber\\
& &\ \ \ \ \cdot \psi_{1}(y_{1})\cdots \psi_{r}(y_{r})
(Y_{\cal{E}}(a,x_{0})b)(x_{2})\phi_{1}(z_{1})\cdots \phi_{s}(z_{s})
\nonumber\\
&=&\Res_{x_{1}}\Res_{x_{0}}x_{0}^{n}\prod_{i=1}^{r}f(x_{2},y_{i})^{k}
\prod_{j=1}^{s}f(x_{2},z_{j})^{k}\nonumber\\
& &\ \ \ \ \cdot \iota_{x_{2},\infty;x_{0},0}(f(x_{2}+x_{0},x_{2})^{-1})
x_{1}^{-1}\delta\left(\frac{x_{2}+x_{0}}{x_{1}}\right)\nonumber\\
& &\ \ \ \ \cdot \left(f(x_{1},x_{2})\psi_{1}(y_{1})\cdots \psi_{r}(y_{r})
a(x_{1})b(x_{2})\phi_{1}(z_{1})\cdots
\phi_{s}(z_{s})\right)\nonumber\\
&=&\Res_{x_{1}}\Res_{x_{0}}x_{0}^{n}\prod_{i=1}^{r}f(x_{1}-x_{0},y_{i})^{k}
\prod_{j=1}^{s}f(x_{1}-x_{0},z_{j})^{k}\nonumber\\
& &\ \ \ \ \cdot \iota_{x_{2},\infty;x_{0},0}(f(x_{2}+x_{0},x_{2})^{-1})
x_{1}^{-1}\delta\left(\frac{x_{2}+x_{0}}{x_{1}}\right)\nonumber\\
& &\ \ \ \ \cdot \left(f(x_{1},x_{2})\psi_{1}(y_{1})\cdots \psi_{r}(y_{r})
a(x_{1})b(x_{2})\phi_{1}(z_{1})\cdots
\phi_{s}(z_{s})\right)\nonumber\\
&=&\Res_{x_{1}}\Res_{x_{0}}x_{0}^{n}
e^{-x_{0}\frac{\partial}{\partial x_{1}}}
\left(\prod_{i=1}^{r}f(x_{1},y_{i})
\prod_{j=1}^{s}f(x_{1},z_{j})\right)^{k}\nonumber\\
& &\ \ \ \ \cdot \iota_{x_{2},\infty;x_{0},0}(f(x_{2}+x_{0},x_{2})^{-1})
x_{1}^{-1}\delta\left(\frac{x_{2}+x_{0}}{x_{1}}\right)\nonumber\\
& &\ \ \ \ \cdot \left(f(x_{1},x_{2})\psi_{1}(y_{1})\cdots \psi_{r}(y_{r})
a(x_{1})b(x_{2})\phi_{1}(z_{1})\cdots
\phi_{s}(z_{s})\right)\nonumber\\
&=&\Res_{x_{1}}\Res_{x_{0}}\sum_{t=0}^{k-1}\frac{(-1)^{t}}{t!}x_{0}^{n+t}
\left(\frac{\partial}{\partial x_{1}}\right)^{t}
\left(\prod_{i=1}^{r}f(x_{1},y_{i})
\prod_{j=1}^{s}f(x_{1},z_{j})\right)^{k}\nonumber\\
& &\ \ \ \ \cdot \iota_{x_{2},\infty;x_{0},0}(f(x_{2}+x_{0},x_{2})^{-1})
x_{1}^{-1}\delta\left(\frac{x_{2}+x_{0}}{x_{1}}\right)\nonumber\\
& &\ \ \ \ \cdot \left(f(x_{1},x_{2})\psi_{1}(y_{1})\cdots \psi_{r}(y_{r})
a(x_{1})b(x_{2})\phi_{1}(z_{1})\cdots
\phi_{s}(z_{s})\right).
\end{eqnarray}
Notice that for any polynomial $B$ and for $0\le t\le k-1$,
$\left(\frac{\partial}{\partial x_{1}}\right)^{t}B^{k}$
is a multiple of $B$.
Using (\ref{elong-exp}) we have
\begin{eqnarray*}
& &P Q R \prod_{i=1}^{r}f(x_{2},y_{i})
\prod_{j=1}^{s}f(x_{2},z_{j})
\sum_{t=0}^{k-1}\frac{(-1)^{t}}{t!}x_{0}^{n+t}
\left(\frac{\partial}{\partial x_{1}}\right)^{t}
\left(\prod_{i=1}^{r}f(x_{1},y_{i})
\prod_{j=1}^{s}f(x_{1},z_{j})\right)^{k}\nonumber\\
& &\ \ \ \ \cdot \iota_{x_{2},\infty;x_{0},0}(f(x_{2}+x_{0},x_{2})^{-1})
x_{1}^{-1}\delta\left(\frac{x_{2}+x_{0}}{x_{1}}\right)\nonumber\\
& &\ \ \ \ \cdot \left(f(x_{1},x_{2})\psi_{1}(y_{1})\cdots \psi_{r}(y_{r})
a(x_{1})b(x_{2})\phi_{1}(z_{1})\cdots
\phi_{s}(z_{s})\right)\nonumber\\
&\in& \left(\Hom (W,W((y_{1}^{-1},\dots,
y_{r}^{-1},x_{2}^{-1},z_{1}^{-1},\dots,z_{s}^{-1})))\right)
((x_{0}))[[x_{1}^{\pm 1}]].
\end{eqnarray*}
Then
\begin{eqnarray}
& &P Q R\prod_{i=1}^{r}f(x_{2},y_{i})^{k+1}
\prod_{j=1}^{s}f(x_{2},z_{j})^{k+1}
\psi_{1}(y_{1})\cdots\psi_{r}(y_{r})
(a(x)_{n}b(x))(x_{2})\phi_{1}(z_{1})\cdots\phi_{s}(z_{s})\nonumber\\
&\in& \Hom (W,W((y_{1}^{-1},\dots,
y_{r}^{-1},x_{2}^{-1},z_{1}^{-1},\dots,z_{s}^{-1}))).
\end{eqnarray}
This proves that $(\psi_{1}(x),\dots,\psi_{r}(x),a(x)_{n}b(x),
\phi_{1}(x), \dots,\phi_{s}(x))$ is quasi-compatible.
It is clear that the assertion with compatibility holds.
\end{proof}

The following is our second key result:

\bt{tmaximal}
Every maximal quasi-compatible subspace of $\E^{o}(W)$
is closed and contains $1_{W}$. Furthermore, for any
quasi-compatible subset $S$, there exists a (unique) smallest closed
quasi-compatible subspace $\<S\>$ containing $S$ and $1_{W}$, and
$(\<S\>,Y_{\E^{o}},1_{W})$ carries the structure of a nonlocal
vertex algebra with $W$ as a faithful (left) quasi module-at-infinity
where the vertex operator map $Y_{W}$ is given by
$Y_{W}(\psi(x),x_{0})=\psi(x_{0})$. If $S$ is compatible, then
$W$ is a module-at-infinity for $\<S\>$.
\et

\begin{proof} Let $K$ be any maximal quasi-compatible subspace of $\E^{o}(W)$.
Clearly, $K+\C 1_{W}$ is quasi-compatible. With $K$ maximal we must
have $1_{W}\in K$. Let $a(x),b(x)\in K,\;n\in \Z$.
It follows from Proposition \ref{pgeneratingcomplicatedone}
and an induction that any finite sequence in $K\cap \{ a(x)_{n}b(x)\}$
is  quasi-compatible. Again, with $K$ maximal we must have
$a(x)_{n}b(x)\in K$. This proves that $K$ is closed.
The rest assertions follow immediately from Theorem
\ref{tclosed}.
\end{proof}

Recall that $\C(x)$ denotes the field of rational functions and
$\iota_{x,0}$ and $\iota_{x,\infty}$ are the field embeddings of
$\C(x)$ into $\C((x))$ and $\C((x^{-1}))$, respectively.

\bp{pconvert}
Let $V$ be a nonlocal vertex algebra generated by
a quasi-compatible subset of $\E^{o}(W)$.
Suppose that the following relation holds
\begin{eqnarray}
& &(x_{1}-x_{2})^{k}p(x_{1},x_{2})a(x_{1})b(x_{2})\nonumber\\
&=&(x_{1}-x_{2})^{k}p(x_{1},x_{2})
\sum_{i=1}^{r}\iota_{x,\infty}(q_{i})(x_{1}-x_{2})
u_{i}(x_{2})v_{i}(x_{1}),
\end{eqnarray}
where $a(x),b(x),u_{i}(x),v_{i}(x)\in V$ and
$p(x,y)\in \C[x,y],\; q_{i}(x)\in \C(x),\; k\in \N$ with $p(x,x)\ne 0$.
Then there exists a nonnegative integer $k'$ such that
\begin{eqnarray}\label{epqsum}
& &(x_{1}-x_{2})^{k'}Y_{\E^{o}}(a(x),x_{1})Y_{\E^{o}}(b(x),x_{2})
\nonumber\\
&=&(x_{1}-x_{2})^{k'}\sum_{i=1}^{r}\iota_{x,0}(q_{i})(-x_{2}+x_{1})
Y_{\E^{o}}(u_{i}(x),x_{2})Y_{\E^{o}}(v_{i}(x),x_{1}).
\end{eqnarray}
\ep

\begin{proof} Let $\theta(x)\in V$. By Lemma \ref{lclosed},
there exists a nonzero polynomial $f(x,y)$ such that
\begin{eqnarray}
& &f(x+x_{1},x)f(x+x_{1},x+x_{2})f(x+x_{2},x)
Y_{\E^{o}}(a(x),x_{1})Y_{\E^{o}}(b(x),x_{2})\theta (x)\nonumber\\
&=&\left(f(y,x)f(y,z)f(z,x)a(y)b(z)\theta (x)\right)|_{y=x+x_{1},z=x+x_{2}}
\end{eqnarray}
and such that
\begin{eqnarray}
& &f(x+x_{1},x)f(x+x_{1},x+x_{2})f(x+x_{2},x)
Y_{\E^{o}}(u_{i}(x),x_{1})Y_{\E^{o}}(v_{i}(x),x_{2})\theta (x)\nonumber\\
&=&\left(f(y,x)f(y,z) f(z,x)
u_{i}(y)v_{i}(z)\theta (x)\right)|_{y=x+x_{1},z=x+x_{2}}
\end{eqnarray}
for $i=1,\dots,r$. Let $0\ne g(x)\in \C[x]$ such that
$g(x)q_{i}(x)\in \C[x]$ for $i=1,\dots,r$.
Then
\begin{eqnarray}
& &f(x+x_{1},x)f(x+x_{1},x+x_{2})f(x+x_{2},x)
(x_{1}-x_{2})^{k}p(x+x_{1},x+x_{2})\nonumber\\
& &\ \ \ \ \cdot g(x_{1}-x_{2})
Y_{\E^{o}}(a(x),x_{1})Y_{\E^{o}}(b(x),x_{2})\theta (x)\nonumber\\
&=&(x_{1}-x_{2})^{k}
\left(f(y,x)f(y,z)f(z,x)g(y-z)p(y,z)a(y)b(z)\theta
(x)\right)|_{y=x+x_{1},z=x+x_{2}}\nonumber\\
&=&(x_{1}-x_{2})^{k}
\left(f(y,x)f(y,z)f(z,x)p(y,z)\sum_{i=1}^{r}(gq_{i})(y-z)
u_{i}(y)v_{i}(z)\theta (x)\right)|_{y=x+x_{1},z=x+x_{2}}
\nonumber\\
&=&f(x+x_{1},x)f(x+x_{1},x+x_{2})f(x+x_{2},x)
(x_{1}-x_{2})^{k}p(x+x_{1},x+x_{2})\nonumber\\
& &\ \ \ \ \cdot\sum_{i=1}^{r}(gq_{i})(x_{1}-x_{2})
Y_{\E^{o}}(u_{i}(x),x_{2})Y_{\E^{o}}(v_{i}(x),x_{1})\theta
(x)\nonumber\\
&=&f(x+x_{1},x)f(x+x_{1},x+x_{2})f(x+x_{2},x)
(x_{1}-x_{2})^{k}p(x+x_{1},x+x_{2})\nonumber\\
& &\ \ \ \ g(x_{1}-x_{2})
\cdot\sum_{i=1}^{r}\iota_{x,0}(q_{i})(-x_{2}+x_{1}))
Y_{\E^{o}}(u_{i}(x),x_{2})Y_{\E^{o}}(v_{i}(x),x_{1})\theta.
\end{eqnarray}
Notice that we can multiply both sides by
$\iota_{x,\infty;x_{1},0}f(x+x_{1},x)^{-1}
\iota_{x,\infty;x_{2},0}f(x+x_{2},x)^{-1}$ to cancel the factors
$f(x+x_{1},x)$ and $f(x+x_{2},x)$ (recall Remark
\ref{rcancellation}). Since $p(x,x)\ne 0$, we can also cancel the
factor $p(x+x_{1},x+x_{2})$. By cancelation we get
\begin{eqnarray*}
& &(x_{1}-x_{2})^{k}f(x+x_{1},x+x_{2})g(x_{1}-x_{2})
Y_{\E^{o}}(a(x),x_{1})Y_{\E^{o}}(b(x),x_{2})\theta (x)\nonumber\\
&=&(x_{1}-x_{2})^{k}f(x+x_{1},x+x_{2})g(x_{1}-x_{2})\\
& &\ \ \cdot \sum_{i=1}^{r}\iota_{x,0}(q_{i})(-x_{2}+x_{1}))
Y_{\E^{o}}(u_{i}(x),x_{2})Y_{\E^{o}}(v_{i}(x),x_{1})\theta (x).\ \ \
\
\end{eqnarray*}
Write $g(x)=x^{l}\bar{g}(x)$, where $l\ge 0,\; \bar{g}(x)\in \C[x]$
with $\bar{g}(0)\ne 0$. Similarly, write
$f(x,z)=(x-z)^{t}\bar{f}(x,z)$, where $t\ge 0,\; \bar{f}(x,z)\in
\C[x,z]$ with $\bar{f}(x,x)\ne 0$.
By a further cancelation we get (\ref{epqsum}) with $k'=k+l+t$.
\end{proof}

Recall that $G$ denotes the group of linear transformations on $\C$.

\bl{lconnection} Let $\Gamma$ be a group of linear transformations
and let $V$ be a vertex algebra generated by a quasi-compatible
subset of $\E^{o}(W)$. Assume that
$$R_{g}a(x)\;(=a(g(x)))\in V\;\;\;\mbox{ for } g\in \Gamma,\; a(x)\in V.$$
Then
\begin{eqnarray}\label{einvariance}
Y_{\E^{o}}(R_{g}a(x),x_{0})R_{g}b(x)=R_{g}Y_{\E^{o}}(a(x),g_{0}x_{0})b(x)
\end{eqnarray}
for $g\in \Gamma,\; a(x),b(x)\in V$, where $g(x)=g_{0}x+g_{1}$. \el

\begin{proof} Let $g\in \Gamma,\;a(x),b(x)\in V$.
There exists a nonzero polynomial $p(x_{1},x_{2})$
such that
$$p(x_{1},x_{2})a(x_{1})b(x_{2})\in
\Hom(W,W((x_{1}^{-1},x_{2}^{-1}))).$$ Then
\begin{eqnarray*}
p(x+g_{0}x_{0},x)Y_{\E^{o}}(a(x),g_{0}x_{0})b(x)
=\left(p(x_{1},x)a(x_{1})b(x)\right)|_{x_{1}=x+g_{0}x_{0}}.
\end{eqnarray*}
Substituting $x$ with $g(x)$ we get
\begin{eqnarray*}
&&p(g(x)+g_{0}x_{0},g(x))R_{g(x)}
\left(Y_{\E^{o}}(a(x),g_{0}x_{0})b(x)\right)\\
 & =&
\left(p(x_{1},g(x))a(x_{1})b(g(x))
\right)|_{x_{1}=g(x)+g_{0}x_{0}=g(x+x_{0})}\\
&=&\left(p(g(x_{1}),g(x))a(g(x_{1})b(x)\right)|_{x_{1}=x+x_{0}}.
\end{eqnarray*}
 We also have
$$p(g(x_{1}),g(x_{2}))a(g(x_{1}))b(g(x_{2}))\in
\Hom(W,W((x_{1}^{-1},x_{2}^{-1}))),$$ so that
$$p(g(x)+g_{0}x_{0},g(x))Y_{\E^{o}}(a(g(x)),x_{0})b(g(x))
=\left(p(g(x_{1}),g(x))a(g(x_{1}))b(g(x))\right)|_{x_{1}=x+x_{0}}.$$
Consequently,
$$p(g(x)+g_{0}x_{0},g(x))Y_{\E^{o}}(a(g(x)),x_{0})b(g(x))
=p(g(x)+g_{0}x_{0},g(x))R_{g(x)}Y_{\E}^{o}(a(x),g_{0}x_{0})b(x).$$
By cancelation, we obtain (\ref{einvariance}).
\end{proof}

The following is an analogue of (\cite{li-gamma}, Proposition 4.3):

\bp{pgamma-locality-key}
Let $a(x),b(x),c(x)\in \E^{o}(W)$. Assume that
\begin{eqnarray*}
&&f(x_{1},x_{2})a(x_{1})b(x_{2})=f(x_{1},x_{2})b(x_{2})a(x_{1}),
\label{efab}\\
&&g(x_{1},x_{2})a(x_{1})c(x_{2})=\tilde{g}(x_{1},x_{2})c(x_{2})a(x_{1}),
\label{egca}\\
&&h(x_{1},x_{2})b(x_{1})c(x_{2})=\tilde{h}(x_{1},x_{2})c(x_{2})b(x_{1}),
\label{ehcb}
\end{eqnarray*}
where $f(x,y), g(x,y),\tilde{g}(x,y), h(x,y), \tilde{h}(x,y)$
 are nonzero polynomials.
Then for any $n\in \Z$,
there exists $k\in \N$, depending on $n$, such that
\begin{eqnarray}\label{einduction-locality}
f(x_{3},x)^{k}g(x_{3},x)a(x_{3})(b(x)_{n}c(x))
=f(x_{3},x)^{k}\tilde{g}(x_{3},x)(b(x)_{n}c(x))a(x_{3}).
\end{eqnarray}
\ep

\begin{proof} Let $n\in \Z$ be arbitrarily fixed.
Let $k$ be a nonnegative integer such that
$$x_{0}^{k+n}\iota_{x,\infty;x_{0},0}(h(x+x_{0},x)^{-1})
\in \C((x^{-1}))[[x_{0}]].$$
In the proof of Proposition 4.3 of \cite{li-gamma}, take
$\alpha=1$ and replace $\iota_{x,x_{0}}$ with $\iota_{x,\infty'x_{0},0}$.
Then the same arguments prove (\ref{einduction-locality}).
\end{proof}

\section{Associative algebra $DY_{q}^{\infty}(sl_{2})$ and
modules-at-infinity for quantum vertex algebras}

In this section we continue to study $\S$-local subsets of
$\E^{o}(W)$ for a vector space $W$ and we prove that any $\S$-local
subset generates a weak quantum vertex algebra with $W$ as a
canonical module-at-infinity. We introduce another version
$DY_{q}^{\infty}(sl_{2})$ of the double Yangian and we prove that
every $DY_{q}^{\infty}(sl_{2})$-module $W$ is naturally a
module-at-infinity for the quantum vertex algebra $V_{q}$ which was
constructed in Section 2.

First we prove a simple result that we shall need later:

\bl{lD-infinity}
Let $V$ be a nonlocal vertex algebra and let $(W,Y_{W})$
be a (left) quasi $V$-module-at-infinity. Then
\begin{eqnarray}\label{e3.1}
Y_{W}(\D v,x)=\frac{d}{dx}Y_{W}(v,x)\ \ \ \mbox{ for }v\in V.
\end{eqnarray}
\el

\begin{proof} For any $v\in V$, by definition there exists a
nonzero polynomial $p(x_{1},x_{2})$ such that
\begin{eqnarray*}
& &p(x_{1},x_{2})Y_{W}(v,x_{1})Y_{W}({\bf 1},x_{2})
\in \Hom (W,W((x_{1}^{-1},x_{2}^{-1}))),\\
& & p(x_{2}+x_{0},x_{2})Y_{W}(Y(v,x_{0}){\bf 1},x_{2})
=\left( p(x_{1},x_{2})Y_{W}(v,x_{1})Y_{W}({\bf
  1},x_{2})\right)|_{x_{1}=x_{2}+x_{0}}.
\end{eqnarray*}
With $Y(v,x_{0}){\bf 1}=e^{x_{0}\D}v$ and $Y_{W}({\bf
1},x_{2})=1_{W}$, we get
\begin{eqnarray*}
p(x_{2}+x_{0},x_{2})Y_{W}(e^{x_{0}\D}v,x_{2})
&=&\left( p(x_{1},x_{2})Y_{W}(v,x_{1})\right)|_{x_{1}=x_{2}+x_{0}}\\
&=&p(x_{2}+x_{0},x_{2})Y_{W}(v,x_{2}+x_{0}).
\end{eqnarray*}
As both $Y_{W}(e^{x_{0}\D}v,x_{2})$ and $Y_{W}(v,x_{2}+x_{0})$ lie
in $(\Hom (W,W((x_{2}^{-1}))))[[x_{0}]]$, in view of Remark
\ref{rcancellation} we have
\begin{eqnarray*}
Y_{W}(e^{x_{0}\D}v,x_{2})=Y_{W}(v,x_{2}+x_{0})
=e^{x_{0}\frac{d}{dx_{2}}}Y_{W}(v,x_{2}),
\end{eqnarray*}
which implies (\ref{e3.1}).
\end{proof}

The following is straightforward to prove:

\bl{lsimple-factx} Let $W$ be any vector space and let
\begin{eqnarray*}
 & &A(x_{1},x_{2})\in \Hom (W,W((x_{2}^{-1}))((x_{1}^{-1}))),\ \
B(x_{1},x_{2})\in \Hom (W,W((x_{1}^{-1}))((x_{2}^{-1}))),\\
& & \ \ \ \ \ \ \ C(x_{2},x_{0})\in (\Hom
(W,W((x_{2}^{-1})))((x_{0})).
\end{eqnarray*}
Then
\begin{eqnarray}
& &x_{0}^{-1}\delta\left(\frac{x_{1}-x_{2}}{x_{0}}\right)
A(x_{1},x_{2})-x_{0}^{-1}\delta\left(\frac{x_{2}-x_{1}}{-x_{0}}\right)
B(x_{1},x_{2})\nonumber\\
& &\ \ \ =x_{2}^{-1}\delta\left(\frac{x_{1}-x_{0}}{x_{2}}\right)
C(x_{2},x_{0})
\end{eqnarray}
if and only if there exist a nonnegative integer $k$ and
$$F(x_{1},x_{2})\in \Hom(W,W((x_{1}^{-1},x_{2}^{-1})))$$
such that
\begin{eqnarray}
& &(x_{1}-x_{2})^{k}A(x_{1},x_{2})=F(x_{1},x_{2})=(x_{1}-x_{2})^{k}B(x_{1},x_{2}),\\
& &x_{0}^{k}C(x_{2},x_{0})=F(x_{2}+x_{0},x_{2}).
\end{eqnarray}
\el

\br{rrestrictions} {\em Let $g(x)=\sum_{n\ge m}g_{n}x^{-n}\in
\C((x^{-1}))$ with $m\in \Z$. We have
$$g(x_{1}-x_{2})=\sum_{n\ge m}g_{n}(x_{1}-x_{2})^{-n}
=\sum_{n\ge m}\sum_{i\in
\N}\binom{-n}{i}(-1)^{i}g_{n}x_{1}^{-n-i}x_{2}^{i} \in
\C[x_{2}]((x_{1}^{-1})).$$ Furthermore, for any $\psi(x),\phi(x)\in
\Hom (W,W((x^{-1})))$ with $W$ a vector space, the product
$$g(x_{1}-x_{2})\psi(x_{2})\phi(x_{1})\ \ \mbox{ exists in }
\Hom (W,W((x_{2}^{-1}))((x_{1}^{-1})))$$
 and we have
\begin{eqnarray*}
x_{0}^{-1}\delta\left(\frac{x_{1}-x_{2}}{-x_{0}}\right)
g(-x_{0})\psi(x_{2})\phi(x_{1})
=x_{0}^{-1}\delta\left(\frac{x_{1}-x_{2}}{-x_{0}}\right)
g(x_{1}-x_{2})\psi(x_{2})\phi(x_{1}).
\end{eqnarray*}}
\er

Using Lemma \ref{lsimple-factx} we immediately have:

\bl{lopp-jacobi}
Let $V$ be a nonlocal vertex algebra, let $(W,Y_{W})$ be a
$V$-module-at-infinity, and let
$$u,v, u^{(i)},v^{(i)}\in V,\; f_{i}(x)\in \C(x)
\ \ (i=1,\dots,r).$$
Then
\begin{eqnarray*}
& &-x_{0}^{-1}\delta\left(\frac{x_{2}-x_{1}}{-x_{0}}\right)
Y_{W}(u,x_{1})Y_{W}(v,x_{2})\nonumber\\
& &\hspace{1cm}
+x_{0}^{-1}\delta\left(\frac{x_{1}-x_{2}}{x_{0}}\right)
\sum_{i=1}^{r}\iota_{x,\infty}(f_{i})(-x_{0})
Y_{W}(v_{(i)},x_{2})Y_{W}(u_{(i)},x_{1})\nonumber\\
&=&x_{2}^{-1}\delta\left(\frac{x_{1}-x_{0}}{x_{2}}\right)
Y_{W}(Y(u,x_{0})v,x_{2})
\end{eqnarray*}
if and only if there exists a nonnegative integer $k$ such that
\begin{eqnarray*}
& &(x_{1}-x_{2})^{k}Y_{W}(u,x_{1})Y_{W}(v,x_{2})\\
&=&(x_{1}-x_{2})^{k}
\sum_{i=1}^{r}\iota_{x,\infty}(f_{i})(-x_{1}+x_{2})
Y_{W}(v^{(i)},x_{2})Y_{W}(u^{(i)},x_{1}).
\end{eqnarray*}
 \el

The following is an analogue of (\cite{li-qva1}, Proposition 6.7):

\bp{pqva-module-infty} Let $V$ be a nonlocal vertex algebra, let
$(W,Y_{W})$ be a $V$-module-at-infinity, and let $$n\in \Z,\; u,v,
u^{(i)},v^{(i)}\in V,\; f_{i}(x)\in \C(x) \ \ (i=1,\dots,r), \;
c^{(0)},\dots,c^{(s)}\in V.$$
 If
\begin{eqnarray}\label{ecross-va}
& &(x_{1}-x_{2})^{n}Y(u,x_{1})Y(v,x_{2})-(-x_{2}+x_{1})^{n}
\sum_{i=1}^{r}\iota_{x,0}(f_{i})(x_{2}-x_{1})
Y(v^{(i)},x_{2})Y(u^{(i)},x_{1})\nonumber\\
&=&\sum_{j=0}^{s} Y(c^{(j)},x_{2})
\frac{1}{j!}\left(\frac{\partial}{\partial x_{2}}\right)^{j}
x_{1}^{-1}\delta\left(\frac{x_{2}}{x_{1}}\right)
\end{eqnarray}
on $V$, then
\begin{eqnarray}\label{ecross-module}
& &(x_{1}-x_{2})^{n}
\sum_{i=1}^{r}\iota_{x,\infty}(f_{i})(-x_{1}+x_{2})
Y_{W}(v^{(i)},x_{2})Y_{W}(u^{(i)},x_{1})\nonumber\\
& &\ \ \ \ \ \ \ -(-x_{2}+x_{1})^{n}Y_{W}(u,x_{1})Y_{W}(v,x_{2})\nonumber\\
&=&\sum_{j=0}^{s} Y_{W}(c^{(j)},x_{2})
\frac{1}{j!}\left(\frac{\partial}{\partial x_{2}}\right)^{j}
x_{1}^{-1}\delta\left(\frac{x_{2}}{x_{1}}\right)
\end{eqnarray}
on $W$. If $(W,Y_{W})$ is faithful, the converse is also true.
\ep

\begin{proof}
Let $k$ be a nonnegative integer such that $k>s$ and
$n+k\ge 0$. {}From (\ref{ecross-va}) we get
$$(x_{1}-x_{2})^{k+n}Y(u,x_{1})Y(v,x_{2})=\sum_{i=1}^{r}(x_{1}-x_{2})^{k+n}
\iota_{x,0}(f_{i})(x_{2}-x_{1})Y(v^{(i)},x_{2})Y(u^{(i)},x_{1}).$$
By Corollary 5.3 in \cite{li-qva2} (cf. \cite{ek}) we have
the $\S$-skew symmetry
\begin{eqnarray}
Y(u,x)v=\sum_{i=1}^{r}\iota_{x,0}(f_{i}(-x))e^{x\D}Y(v^{(i)},-x)u^{(i)}.
\end{eqnarray}
We also have the following $\S$-Jacobi identity
\begin{eqnarray}\label{e4.10proof}
& &x_{0}^{-1}\delta\left(\frac{x_{1}-x_{2}}{x_{0}}\right)
Y(u,x_{1})Y(v,x_{2})\nonumber\\
& &\hspace{1cm}
-x_{0}^{-1}\delta\left(\frac{x_{2}-x_{1}}{-x_{0}}\right)
\sum_{i=1}^{r}\iota_{x_{0},0}(f_{i}(-x_{0}))
Y(v^{(i)},x_{2})Y(u^{(i)},x_{1})\nonumber\\
&=&x_{2}^{-1}\delta\left(\frac{x_{1}-x_{0}}{x_{2}}\right)
Y(Y(u,x_{0})v,x_{2}).
\end{eqnarray}
By taking $\Res_{x_{0}}x_{0}^{n}$ we get
\begin{eqnarray}
& &(x_{1}-x_{2})^{n}Y(u,x_{1})Y(v,x_{2})-(-x_{2}+x_{1})^{n}
\sum_{i=1}^{r}\iota_{x,0}(f_{i})(x_{2}-x_{1})
Y(v^{(i)},x_{2})Y(u^{(i)},x_{1})\ \ \ \nonumber\\
& &\ \ \ \ =\sum_{j\ge 0}Y(u_{n+j}v,x_{2})
\frac{1}{j!}\left(\frac{\partial}{\partial x_{2}}\right)^{j}
x_{1}^{-1}\delta\left(\frac{x_{2}}{x_{1}}\right).
\end{eqnarray}
Combining this with (\ref{ecross-va}) we obtain
\begin{eqnarray}\label{euv=c}
u_{n+j}v=c^{(j)}\;\;\;\mbox{ for }j=0,\dots,s,\;\;\mbox{ and }\;
u_{n+j}v=0\;\;\;\mbox{ for }j>s.
\end{eqnarray}

Let $l$ be a sufficiently large nonnegative integer such that
\begin{eqnarray*}
& &(x_{1}-x_{2})^{l}Y_{W}(u,x_{1})Y_{W}(v,x_{2})=F(x_{1},x_{2})
\in \Hom (W,W((x_{1}^{-1},x_{2}^{-1}))),\\
& &x_{0}^{l}Y_{W}(Y(u,x_{0})v,x_{2})
=F(x_{2}+x_{0},x_{2})
\end{eqnarray*}
and
\begin{eqnarray*}
& &(x_{1}-x_{2})^{l}Y_{W}(v^{(i)},x_{2})Y_{W}(u^{(i)},x_{1})=G_{i}(x_{1},x_{2})
\in \Hom (W,W((x_{1}^{-1},x_{2}^{-1}))),\\
& &x_{0}^{l}Y_{W}(Y(v^{(i)},-x_{0})u^{(i)},x_{1})
=G_{i}(x_{1},x_{1}-x_{0})
\end{eqnarray*}
for all $i=1,\dots,r$.

Let $p(x)$ be a nonzero polynomial such that
$$p(x)f_{i}(x)\in \C[x]\ \ \ \mbox{ for all }i=1,\dots,r.$$

Note that by Lemma \ref{lD-infinity}, we have $Y_{W}(\D
v,x)=(d/dx)Y_{W}(v,x)$ for $v\in V$. Using all of these and the
$\S$-skew symmetry we get
\begin{eqnarray*} p(-x_{0})F(x_{2}+x_{0},x_{2})
&=&p(-x_{0})x_{0}^{l}Y_{W}(Y(u,x_{0})v,x_{2})\\
&=&\sum_{i=1}^{r}p(-x_{0})f_{i}(-x_{0})x_{0}^{l}
Y_{W}(e^{x_{0}\D}Y(v^{(i)},-x_{0})u^{(i)},x_{2})\\
&=&\sum_{i=1}^{r}p(-x_{0})f_{i}(-x_{0})x_{0}^{l}
Y_{W}(Y(v^{(i)},-x_{0})u^{(i)},x_{2}+x_{0})\\
&=&\sum_{i=1}^{r}p(-x_{0})f_{i}(-x_{0})
\left( G_{i}(x_{1},x_{1}-x_{0})\right)|_{x_{1}=x_{2}+x_{0}}\\
&=&\sum_{i=1}^{r}p(-x_{0})f_{i}(-x_{0})
G_{i}(x_{2}+x_{0},x_{2}),
\end{eqnarray*}
which implies
\begin{eqnarray*}
p(x_{2}-x_{1})F(x_{1},x_{2})=\sum_{i=1}^{r}(pf_{i})(x_{2}-x_{1})
G_{i}(x_{1},x_{2}).
\end{eqnarray*}
Then
\begin{eqnarray*}
& &p(x_{2}-x_{1})(x_{1}-x_{2})^{l}Y_{W}(u,x_{1})Y_{W}(v,x_{2})\\
&=&p(x_{2}-x_{1})
\sum_{i=1}^{r}(x_{1}-x_{2})^{l}\iota_{x,\infty}(f_{i})(-x_{1}+x_{2})
Y_{W}(v^{(i)},x_{2})Y_{W}(u^{(i)},x_{1}).
\end{eqnarray*}
As both
$(x_{1}-x_{2})^{l}Y_{W}(u,x_{1})Y_{W}(v,x_{2})$ and
$$\sum_{i=1}^{r}(x_{1}-x_{2})^{l}\iota_{x,\infty}(f_{i})(-x_{1}+x_{2})
Y_{W}(v^{(i)},x_{2})Y_{W}(u^{(i)},x_{1})$$ lie in
$\Hom(W,W((x_{2}^{-1}))((x_{1}^{-1})))$, by Remark
\ref{rcancellation} we get
\begin{eqnarray*}
& &(x_{1}-x_{2})^{l}Y_{W}(u,x_{1})Y_{W}(v,x_{2})\\
&=&(x_{1}-x_{2})^{l}\sum_{i=1}^{r}\iota_{x,\infty}(f_{i})(-x_{1}+x_{2})
Y_{W}(v^{(i)},x_{2})Y_{W}(u^{(i)},x_{1}).
\end{eqnarray*}
Now by Lemma \ref{lsimple-factx} we have
\begin{eqnarray}\label{eS-jacobi-module}
& &-x_{0}^{-1}\delta\left(\frac{x_{2}-x_{1}}{-x_{0}}\right)
Y_{W}(u,x_{1})Y_{W}(v,x_{2})\nonumber\\
& &\hspace{1cm}
+x_{0}^{-1}\delta\left(\frac{x_{1}-x_{2}}{x_{0}}\right)
\sum_{i=1}^{r}\iota_{x_{0},\infty}(f_{i}(-x_{0}))
Y_{W}(v^{(i)},x_{2})Y_{W}(u^{(i)},x_{1})\nonumber\\
&=&x_{2}^{-1}\delta\left(\frac{x_{1}-x_{0}}{x_{2}}\right)
Y_{W}(Y(u,x_{0})v,x_{2}).
\end{eqnarray}
Using this and (\ref{euv=c}) we obtain (\ref{ecross-module}).

For the converse, we trace back, assuming that $W$ is faithful and
(\ref{ecross-module}) holds. Let $k$ be a nonnegative integer such
that $k+n\ge 0$ and $k>s$. Then
\begin{eqnarray}\label{e4.14proof}
& &(x_{1}-x_{2})^{k+n}Y_{W}(u,x_{1})Y_{W}(v,x_{2})\\
&=&(x_{1}-x_{2})^{k+n}\sum_{i=1}^{r}\iota_{x,\infty}(f_{i})(-x_{1}+x_{2})
Y_{W}(v^{(i)},x_{2})Y_{W}(u^{(i)},x_{1}).
\end{eqnarray}
Using Lemma \ref{lsimple-factx} we get (\ref{eS-jacobi-module}).
Combining (\ref{eS-jacobi-module}) with (\ref{ecross-module}) we
obtain (\ref{euv=c}), using the assumption that $Y_{W}$ is
injective. Using (\ref{e4.14proof}) we also have
\begin{eqnarray*}
& &-x_{0}^{-1}\delta\left(\frac{x_{2}-x_{1}}{-x_{0}}\right)
Y_{W}(u,x_{1})Y_{W}(v,x_{2})\nonumber\\
& &\hspace{1cm}
+x_{0}^{-1}\delta\left(\frac{x_{1}-x_{2}}{x_{0}}\right)
\sum_{i=1}^{r}\iota_{x_{0},\infty}(f_{i}(-x_{0}))
Y_{W}(v^{(i)},x_{2})Y_{W}(u^{(i)},x_{1})\nonumber\\
&=&x_{1}^{-1}\delta\left(\frac{x_{2}+x_{0}}{x_{1}}\right)
\sum_{i=1}^{r}\iota_{x_{0},\infty}(f_{i}(-x_{0}))
Y_{W}(Y(v^{(i)},-x_{0})u^{(i)},x_{1})\\
&=&x_{1}^{-1}\delta\left(\frac{x_{2}+x_{0}}{x_{1}}\right)
\sum_{i=1}^{r}\iota_{x_{0},\infty}(f_{i}(-x_{0}))
Y_{W}(Y(v^{(i)},-x_{0})u^{(i)},x_{2}+x_{0})\\
&=&x_{2}^{-1}\delta\left(\frac{x_{1}-x_{0}}{x_{2}}\right)
\sum_{i=1}^{r}\iota_{x_{0},\infty}(f_{i}(-x_{0}))
Y_{W}(e^{x_{0}\D}Y(v^{(i)},-x_{0})u^{(i)},x_{2}).
\end{eqnarray*}
Combining this with (\ref{eS-jacobi-module}) we get the $\S$-skew
symmetry, with which we obtain the $\S$-Jacobi identity
(\ref{e4.10proof}). Then using (\ref{euv=c}) we obtain
(\ref{ecross-va}).
\end{proof}

\bd{ds-locality}
{\em Let $W$ be a vector space. A subset $U$ of $\E^{o}(W)$ is said to be
{\em $\S$-local} if for any $a(x),b(x)\in U$, there exist
$$u^{(1)}(x), v^{(1)}(x),\dots, u^{(r)}(x),v^{(r)}(x)\in U\;\;
\mbox{ and }\ \ f_{1}(x),\dots,f_{r}(x)\in \C(x)$$ such that
\begin{eqnarray}\label{eslocality-relation}
(x_{1}-x_{2})^{k}a(x_{1})b(x_{2})=(x_{1}-x_{2})^{k}
\sum_{i=1}^{r}\iota_{x,\infty}(f_{i})(-x_{1}+x_{2})
u^{(i)}(x_{2})v^{(i)}(x_{1})
\end{eqnarray}
for some nonnegative integer $k$.}
\ed

\bt{tqva-main}
Let $W$ be a vector space and let $U$ be an
$\S$-local subset of $\E^{o}(W)$. Then $U$ is compatible and the
nonlocal vertex algebra $\<U\>$ generated by $U$ is a weak quantum
vertex algebra with $W$ as a module-at-infinity.
\et

\begin{proof} Notice that the relation (\ref{eslocality-relation}) implies that
$$(x_{1}-x_{2})^{k}a(x_{1})b(x_{2})\in \Hom (W,W((x_{1}^{-1},x_{2}^{-1}))).$$
Thus any ordered pair in $U$ is compatible.
As in the proof of Lemma 3.2 of \cite{li-qva1},
using induction we see that any finite sequence in $U$ is compatible.
That is, $U$ is compatible.
By Theorem \ref{tmaximal}, $U$ generates
a nonlocal vertex algebra $\<U\>$ inside $\E^{o}(W)$ with $W$ as a
module-at-infinity. As $U$ is $\S$-local, {}from Proposition \ref{pconvert},
the vertex operators $Y_{\E^{o}}(a(x),x_{0})$ for $a(x)\in U$ form
an $\S$-local subset of $\<U\>$. Because $U$ generates $\<U\>$,
by Lemma 2.7 of \cite{li-qva2}, $\<U\>$ is a weak quantum vertex algebra.
\end{proof}

\bl{ldef-case-Slocal} Let $a(x),b(x)\in \E^{o}(W)$. Suppose that
there exist
$$u^{(i)}(x),v^{(i)}(x)\in \E^{o}(W),\ \ f_{i}(x)\in \C(x) \ (i=1,\dots,r)$$
such that
\begin{eqnarray}
(x_{1}-x_{2})^{k}a(x_{1})b(x_{2})= (x_{1}-x_{2})^{k} \sum_{i=1}^{r}
\iota_{x,\infty}(f_{i})(x_{1}-x_{2}) u^{(i)}(x_{2})v^{(i)}(x_{1})
\end{eqnarray}
for some nonnegative integer $k$. Then $(a(x),b(x))$ is compatible
and
\begin{eqnarray}
& &Y_{\E^{o}}(a(x),x_{0})b(x)\nonumber\\
&=& \Res_{x_{1}}x_{0}^{-1}\delta\left(\frac{x-x_{1}}{x_{0}}\right)
a(x_{1})b(x) -x_{0}^{-1}\delta\left(\frac{x_{1}-x}{-x_{0}}\right)
\sum_{i=1}^{r}\iota_{x,\infty}(f_{i})(-x_{0})
u^{(i)}(x)v^{(i)}(x_{1}).\ \ \ \
\end{eqnarray}
\el

\begin{proof} We have
$$(x_{1}-x_{2})^{k}a(x_{1})b(x_{2})\in
\Hom(W,W((x_{1}^{-1},x_{2}^{-1}))),$$ so that
$$x_{0}^{k}Y_{\E^{o}(W)}(a(x),x_{0})b(x)
=\left((x_{1}-x_{2})^{k}a(x_{1})b(x_{2})\right)|_{x_{1}=x_{2}+x_{0}}.$$
Then we have the Jacobi identity, then the iterate formula.
\end{proof}

Now we come back to double Yangians. Recall that $T$ is the tensor
algebra over the space $sl_{2}\otimes \C[t,t^{-1}]$ and
$T=\coprod_{n\in \Z}T_{n}$ is $\Z$-graded with $\deg (sl_{2}\otimes
t^{n})=n$ for $n\in \Z$. For $n\in \Z$, set $J[n]=\coprod_{m\le
-n}T_{m}$. This gives a decreasing filtration. Denote by $\tilde{T}$
the completion  of $T$ associated with this filtration.

\bd{dyangian-infity}
{\em Let $q$ be a nonzero complex number as before.
We define $DY_{q}^{\infty}(sl_{2})$
to be the quotient algebra of $\tilde{T}$ modulo the following relations:
\begin{eqnarray*}
& &e(x_{1})e(x_{2})
=\frac{x_{1}-x_{2}+q}{x_{1}-x_{2}-q}
e(x_{2})e(x_{1}),\\
& &f(x_{1})f(x_{2})
=\frac{x_{1}-x_{2}-q}{x_{1}-x_{2}+q}
f(x_{2})f(x_{1}),\\
& &[e(x_{1}),f(x_{2})]
=-x_{1}^{-1}\delta\left(\frac{x_{2}}{x_{1}}\right)h(x_{2}),\\
& &h(x_{1})e(x_{2})
=\frac{x_{1}-x_{2}+q}{x_{1}-x_{2}-q}
e(x_{2})h(x_{1}),\\
& &h(x_{1})h(x_{2})=h(x_{2})h(x_{1}),\\
& &h(x_{1})f(x_{2})
=\frac{x_{1}-x_{2}-q}{x_{1}-x_{2}+q}
f(x_{2})h(x_{1}),
\end{eqnarray*}
where it is understood that
$$(x_{1}-x_{2}\pm q)^{-1}=\sum_{i\in \N}(\pm q)^{i}(x_{1}-x_{2})^{-i-1}
\in \C[[x_{1}^{-1},x_{2}]].$$}
\ed

We define a {\em $DY_{q}^{\infty}(sl_{2})$-module} to be a
$T(sl_{2}\otimes \C[t,t^{-1}])$-module $W$ such that for every $w\in W$,
\begin{eqnarray}
sl(n)w=0\ \ \ \mbox{ for $n$ sufficiently small}
\end{eqnarray}
and such that all the defining relations for $DY_{q}^{\infty}(sl_{2})$ hold.
Then for any $DY_{q}^{\infty}(sl_{2})$-module $W$, the generating functions
$e(x),f(x), h(x)$ are elements of $\E^{o}(W)$.
Recall from Section 2 the quantum vertex algebra $V_{q}$. Then we have:

\bt{tdy-infinity} Let $q$ be any nonzero complex number and let $W$
be any $DY_{q}^{\infty}(sl_{2})$-module. There exists one and only
one structure of a $V_{q}$-module-at-infinity on $W$ with
$$Y_{W}(e,x)=e(x),\ \ Y_{W}(f,x)=f(x),\ \  Y_{W}(h,x)=h(x).$$
\et

\begin{proof} The uniqueness is clear as $e,f,h$ generate $V_{q}$.
The proof for the existence is similar to the proof of
Theorem \ref{tmain-dy-0}.
Set $U=\{ e(x),f(x),h(x)\}\subset \E^{o}(W)$. {}From the defining
relations of $DY_{q}^{\infty}(sl_{2})$, $U$ is an $\S$-local subset.
Then, by Theorem \ref{tqva-main}, $U$ generates a
weak quantum vertex algebra $V_{W}$ with $W$ as a faithful
module-at-infinity where $Y_{W}(a(x),x_{0})=a(x_{0})$ for $a(x)\in V_{W}$.
Using Proposition \ref{pqva-module-infty}, we see that
$V_{W}$ is a $DY_{q}(sl_{2})$-module
 with $e(x_{0}), f(x_{0}),h(x_{0})$ acting as
$Y_{\E^{o}}(e(x),x_{0}), Y_{\E^{o}}(f(x),x_{0}),
 Y_{\E^{o}}(h(x),x_{0})$.
Clearly, $1_{W}$ is a vacuum vector of $V_{W}$ viewed as a
$DY_{q}(sl_{2})$-module.
Then $(V_{W},1_{W})$ is a vacuum $DY_{q}(sl_{2})$-module
with an operator $\D$ such that
$\D (1_{W})=0$ and
$$[\D,u(x)]=\frac{d}{dx} u(x)\ \ \ \mbox{ for }u\in sl_{2}.$$
 By the universal property of $V_{q}$,
there exists a $DY_{q}(sl_{2})$-module homomorphism $\theta$ from
$V_{q}$ to $V_{W}$, sending ${\bf 1}$ to $1_{W}$. Since $sl_{2}$
generates $V_{q}$ as a nonlocal vertex algebra, it follows that
$\theta$ is a homomorphism of nonlocal vertex algebras. Using
$\theta$ we obtain a structure of a $V_{q}$-module-at-infinity on
$W$ with the desired property.
\end{proof}

\section{Quasi modules-at-infinity for vertex algebras}
In this section we study quasi-local subsets of $\E^{o}(W)$ for a
vector space $W$ and we prove that every quasi-local subset
generates a vertex algebra with $W$ as a quasi module-at-infinity.
We give a family of examples related to infinite-dimensional Lie
algebras of a certain type, including the Lie algebra of
pseudo-differential operators on the circle.

First we prove:

\bl{lva-quasi-module} Let $V$ be a vertex algebra and let
$(W,Y_{W})$ be a quasi module-at-infinity for $V$ viewed as a
nonlocal vertex algebra. Then for $u,v\in V$, there exists a nonzero
polynomial $p(x_{1},x_{2})$ such that
\begin{eqnarray}\label{epcomm}
p(x_{1},x_{2})Y_{W}(v,x_{2})Y_{W}(u,x_{1})=
p(x_{1},x_{2})Y_{W}(u,x_{1})Y_{W}(v,x_{2}).
\end{eqnarray}
Furthermore, if $(W,Y_{W})$ is a module-at-infinity, then for
$u,v\in V$,
\begin{eqnarray}\label{eqjacobi-lemma}
& &x_{0}^{-1}\delta\left(\frac{x_{1}-x_{2}}{x_{0}}\right)
Y_{W}(v,x_{2})Y_{W}(u,x_{1})
-x_{0}^{-1}\delta\left(\frac{x_{2}-x_{1}}{-x_{0}}\right)
Y_{W}(u,x_{1})Y_{W}(v,x_{2})
\nonumber\\
& &\ \ \ \ \ \ =x_{2}^{-1}\delta\left(\frac{x_{1}-x_{0}}{x_{2}}\right)
Y_{W}(Y(u,x_{0})v,x_{2}).
\end{eqnarray}
\el

\begin{proof} It basically follows from the arguments
of \cite{ll} (Theorem 3.6.3). Let $u,v\in V$. {}From definition,
there exist
$$F(x_{1},x_{2}), G(x_{1},x_{2})\in \Hom(W,W((x_{1}^{-1},x_{2}^{-1}))),\ \
0\ne p(x_{1},x_{2})\in \C[x_{1},x_{2}]$$
such that
$$p(x_{1},x_{2})Y_{W}(u,x_{1})Y_{W}(v,x_{2})=F(x_{1},x_{2}),\ \
p(x_{1},x_{2})Y_{W}(v,x_{2})Y_{W}(u,x_{1})=G(x_{1},x_{2})$$
and
\begin{eqnarray}
& &p(x_{2}+x_{0},x_{2})Y_{W}(Y(u,x_{0})v,x_{2})
=F(x_{2}+x_{0},x_{2})\label{e3.2}\\
& &p(x_{1},x_{1}-x_{0})Y_{W}(Y(v,-x_{0})u,x_{1})
=G(x_{1},x_{1}-x_{0}).\ \ \ \ \label{e3.3}
\end{eqnarray}
Using the skew symmetry of the vertex algebra $V$ and
Lemma \ref{lD-infinity}, we have
$$Y_{W}(Y(v,-x_{0})u,x_{1})=Y_{W}(e^{-x_{0}\D}Y(u,x_{0})v,x_{1})
=Y_{W}(Y(u,x_{0})v,x_{1}-x_{0}).$$
Now (\ref{e3.3}) is rewritten as
$$p(x_{1},x_{1}-x_{0})Y_{W}(Y(u,x_{0})v,x_{1}-x_{0})
=G(x_{1},x_{1}-x_{0}),$$
which gives
$$p(x_{2}+x_{0},x_{2})Y_{W}(Y(u,x_{0})v,x_{2})=G(x_{2}+x_{0},x_{2}).$$
Combining this with (\ref{e3.2}) we get $F(x_{1},x_{2})=G(x_{1},x_{2})$, proving
(\ref{epcomm}).
For the second assertion, the polynomial $p(x_{1},x_{2})$
in the above argument is of the form
$(x_{1}-x_{2})^{k}$ for $k\in \N$. Then it follows from Lemma \ref{lsimple-factx}.
\end{proof}

\br{rvoa-right-module} {\em Let $V$ be a vertex algebra.  A {\em
right $V$-module} (see \cite{hl}, \cite{li-reg}) is a vector space
$W$ equipped with a linear map
$$Y_{W}: V\rightarrow \Hom (W,W((x^{-1})))\subset (\End W)[[x,x^{-1}]]$$
satisfying the following conditions
\begin{eqnarray}
Y_{W}({\bf 1},x)=1_{W},
\end{eqnarray}
and for $u,v\in V$,
\begin{eqnarray}
& &x_{0}^{-1}\delta\left(\frac{x_{1}-x_{2}}{x_{0}}\right)
Y_{W}(v,x_{2})Y_{W}(u,x_{1})
-x_{0}^{-1}\delta\left(\frac{x_{2}-x_{1}}{-x_{0}}\right)
Y_{W}(u,x_{1})Y_{W}(v,x_{2})\nonumber\\
& &=x_{2}^{-1}\delta\left(\frac{x_{1}-x_{0}}{x_{2}}\right)
Y_{W}(Y(u,x_{0})v,x_{2}).
\end{eqnarray}
In view of Lemma \ref{lva-quasi-module}, this notion of right
$V$-module is equivalent to the notion of (left)
$V$-module-at-infinity. } \er

The following is a counterpart of the notion of quasi-local subset
in \cite{li-gamma}:

\bd{dlocality}
{\em  Let $W$ be a vector space. A subset $S$ of $\E^{o}(W)$ is said
to be {\em quasi-local} if for any $a(x),b(x)\in S$,
there exists a nonzero polynomial $p(x_{1},x_{2})$ such that
$$p(x_{1},x_{2})a(x_{1})b(x_{2})=p(x_{1},x_{2})b(x_{2})a(x_{1}).$$
In case that for any $a(x),b(x)\in S$, there exists a nonnegative
integer $k$ such that
$$(x_{1}-x_{2})^{k}a(x_{1})b(x_{2})=(x_{1}-x_{2})^{k}b(x_{2})a(x_{1}),$$
we say $S$ is {\em local}. }
\ed

Specializing Theorem \ref{tmaximal} we have:

\bt{tquasi-local} Let $W$ be a vector space. Every quasi-local
subset $U$ of $\E^{o}(W)$ is quasi-compatible and the nonlocal
vertex algebra $\<U\>$  generated by $U$ inside $\E^{o}(W)$ is a
vertex algebra with $W$ as a (left) quasi module-at-infinity, where
$Y_{W}(a(x),x_{0})=a(x_{0})$ for $a(x)\in \<U\>$. If $U$ is local,
$W$ is a module-at-infinity. \et

\begin{proof} It is clear that
any quasi-local subset $U$ is quasi-compatible. By Theorem
\ref{tmaximal}, $U$ generates a nonlocal vertex algebra $\<U\>$ with
$W$ as a quasi module-at-infinity. {}From Proposition
\ref{pconvert}, $U$ is a local subspace of the nonlocal vertex
algebra $\<U\>$ in the sense that the adjoint vertex operators
associated to the vectors of $U$ are mutually local. As $U$
generates $\<U\>$ as a nonlocal vertex algebra, from \cite{li-g1}
(Proposition 2.17) $\<U\>$ is a vertex algebra. As locality implies
compatibility, the last assertion follows from Theorem
\ref{tmaximal}.
\end{proof}

The following notion was due to \cite{gkk}:

\bd{dgamma-locality} {\em Let $\Gamma$ be a subgroup of $G$ (the
group of linear transformations on $\C$). A subset $S$ of
$\E^{o}(W)$ is said to be {\em $\Gamma$-local} if for any
$a(x),b(x)\in S$, there exist finitely many $g_{1}(x),\dots,
g_{r}(x)\in \Gamma$ such that
$$(x_{1}-g_{1}(x_{2}))\cdots (x_{1}-g_{r}(x_{2}))[a(x_{1}),b(x_{2})]=0.$$ }
\ed

Recall the following notion from \cite{li-gamodule} (cf.
\cite{li-gamma}):

\bd{dgamma-va} {\em Let $\Gamma$ be an abstract group.  A {\em
$\Gamma$-vertex algebra} is a vertex algebra $V$ equipped with two
group homomorphisms
\begin{eqnarray}
R: \Gamma\rightarrow GL(V),\ \ \ \phi: \Gamma\rightarrow \C^{\times}
\end{eqnarray}
such $R_{g}{\bf 1}={\bf 1}$ and
\begin{eqnarray}
R_{g}Y(u,x)v=Y(R_{g}u,\phi(g)^{-1}x)R_{g}v
\end{eqnarray}
for $g\in \Gamma,\; u,v\in V$. }
\ed

\bex{gammava-example} {\em  Let $V=\coprod_{n\in \Z}V_{(n)}$ be a
$\Z$-graded vertex algebra in the sense that $V$ is a vertex algebra
equipped with a $\Z$-grading such that ${\bf 1}\in V_{(0)}$ and
\begin{eqnarray}
u_{m}V_{(n)}\subset
V_{(k+n-m-1)}\ \ \ \mbox{ for }u\in V_{(k)},\; k,m,n\in \Z.
\end{eqnarray}
Let $\Gamma$ be a group acting on $V$ by automorphisms that preserve
the grading and let $\phi: \Gamma \rightarrow \C^{\times}$ be any
group homomorphism. Then $V$ becomes a $\Gamma$-vertex algebra with
$$R_{g}=\phi(g)^{-L(0)}g\ \ \ \mbox{ for }g\in \Gamma,$$
where $L(0)$ denotes the degree operator on $V$ (see
\cite{li-gamodule}). } \eex

Note that the projection $\frac{d}{dx}: G\rightarrow \C^{\times}$ is
a group homomorphism. Then any group homomorphism $\Phi:
\Gamma\rightarrow G$ gives rise to a group homomorphism
$\Phi_{0}=\frac{d}{dx}\circ \Phi:\Gamma\rightarrow \C^{\times}$.

\bd{dgamma-quasimodule-infinity} {\em  Let $V$ be a $\Gamma$-vertex
algebra. A {\em quasi $V$-module-at-infinity} is a quasi
module-at-infinity $(W,Y_{W})$ for $V$ viewed as a nonlocal vertex
algebra, equipped with a group homomorphism
$$\Phi: \Gamma \rightarrow G,  $$
such that $\phi=\Phi_{0}\;(=\frac{d}{dx}\circ \Phi)$ and
\begin{eqnarray}
Y_{W}(R_{g}v,x)=Y_{W}(v,\Phi(g)(x))\;\;\;\mbox{ for }g\in \Gamma,\;
v\in V,
\end{eqnarray}
and $\{Y_{W}(v,x)\;|\; v\in V\}$ is $\Phi(\Gamma)$-local.} \ed

We shall need the following technical result:

\bl{lva-gamma-quasi-module} Let $V$ be a $\Gamma$-vertex algebra,
let $\Phi: \Gamma\rightarrow G$ be a group homomorphism with
$\Phi_{0}=\phi$, and let $(W,Y_{W})$ be a quasi module-at-infinity
for $V$ viewed as a nonlocal vertex algebra. Assume that $\{
Y_{W}(u,x)\;|\; u\in U\}$ is $\Phi(\Gamma)$-local and
$$Y_{W}(R_{g}u,x)=Y_{W}(u,\Phi(g)(x))\ \ \ \mbox{ for }g\in \Gamma,\; u\in U,$$
where $U$ is a $\Gamma$-submodule and a generating subspace of $V$.
Then $(W,Y_{W})$ is a quasi $V$-module-at-infinity. \el

\begin{proof} First we prove that for $u,v\in V$, if $Y_{W}(u,x)$ and
  $Y_{W}(v,x)$ are quasi compatible, then
\begin{eqnarray}\label{eyqumv}
Y_{W}(u_{n}v,x)=Y_{W}(u,x)_{n}Y_{W}(v,x)\ \ \ \mbox{ for }n\in \Z.
\end{eqnarray}
Let $p(x_{1},x_{2})$ be a nonzero polynomial such that
$$p(x_{1},x_{2})Y_{W}(u,x_{1})Y_{W}(v,x_{2})
\in \Hom (W,W((x_{1}^{-1},x_{2}^{-1}))).$$
We have
$$p(x_{0}+x,x)Y_{\E^{o}}(Y_{W}(u,x),x_{0})Y_{W}(v,x)
=\left( p(x_{1},x)Y_{W}(u,x_{1})Y_{W}(v,x)\right)|_{x_{1}=x+x_{0}}.$$
On the other hand, there exists a nonzero polynomial
$q(x_{1},x_{2})$ such that
$$q(x_{0}+x,x)Y_{W}(Y(u,x_{0})v,x)
=\left( q(x_{1},x)Y_{W}(u,x_{1})Y_{W}(v,x)\right)|_{x_{1}=x+x_{0}}.$$
Then
$$p(x_{0}+x)q(x_{0}+x,x)Y_{W}(Y(u,x_{0})v,x)
=p(x_{0}+x)q(x_{0}+x,x)Y_{\E^{o}}(Y_{W}(u,x),x_{0})Y_{W}(v,x).$$
Consequently, we get
$$Y_{W}(Y(u,x_{0})v,x)
=Y_{\E^{o}}(Y_{W}(u,x),x_{0})Y_{W}(v,x),$$ proving (\ref{eyqumv}).
It follows from Proposition \ref{pgamma-locality-key} and induction
that $\{ Y_{W}(v,x)\;|\; v\in V\}$ is $\Phi(\Gamma)$-local.

Suppose that $Y_{W}(R_{g}u,x)=Y_{W}(u,\Phi(g)(x))$ and
$Y_{W}(R_{g}v,x)=Y_{W}(v,\Phi(g)(x))$ for some $g\in \Gamma,\;
u,v\in V$. By suitably choosing a nonzero polynomial
$p(x_{1},x_{2})$ we have
\begin{eqnarray*}
& &p(\phi(g)^{-1}x_{0}+x,x)Y_{W}(R_{g}Y(u,x_{0})v,x)\\
&=&p(\phi(g)^{-1}x_{0}+x,x)Y_{W}(Y(R_{g}u,\phi(g)^{-1}x_{0})R_{g}v,x)\\
&=&\left(p(x_{1},x)Y_{W}(R_{g}u,x_{1})Y_{W}(R_{g}v,x)\right)|_{x_{1}=
x+\phi(g)^{-1}x_{0}}\\
&=&\left(p(x_{1},x)
Y_{W}(u,\Phi(g)(x_{1}))Y_{W}(v,\Phi(g)(x))\right)|_{x_{1}=
x+\phi(g)^{-1}x_{0}}\\
&=&p(\phi(g)^{-1}x_{0}+x,x)Y_{W}(Y(u,x_{0})v,\Phi(g)(x)),
\end{eqnarray*}
which implies
$$Y_{W}(R_{g}Y(u,x_{0})v,x)=Y_{W}(Y(u,x_{0})v,\Phi(g)(x)),$$
noticing
that $\Phi(g)(x+\phi(g)^{-1}x_{0})=\Phi(x)+x_{0}$. Then it follows
{}from induction that
$$Y_{W}(R_{g}v,x)=Y_{W}(v,\Phi(g)(x))\ \ \ \mbox{ for all }g\in \Gamma,\; v\in
V.$$ Thus $W$  is a quasi $V$-module-at-infinity.
\end{proof}

Now we have:

\bt{tmain-2} Let $W$ be a vector space, let $\Gamma$ be a subgroup
of $G$ (the group of linear transformations on $\C$), and let $S$ be
any $\Gamma$-local subset of $\E^{o}(W)$. Set
$$\Gamma\cdot S={\rm span}\{ R_{g}a(x)=a(g(x))\;|\;g\in \Gamma,\; a(x)\in
S\}.$$ Then $\Gamma \cdot S$ is $\Gamma$-local and $\<\Gamma\cdot
S\>$ is a $\Gamma$-vertex algebra with group homomorphisms
$$R: \Gamma\rightarrow GL(\<\Gamma\cdot S\>) \ \ \mbox{ and }\ \
\phi: \Gamma\rightarrow \C^{\times}$$
 defined by
\begin{eqnarray}
& &R_{g}(\alpha(x))=\alpha (g(x))\;\;\;\mbox{ for }g\in \Gamma,\;
\alpha(x)\in \<\Gamma\cdot S\>,\\
& &\phi(g(x))=g_{0}\;\;\;\mbox{ for }g(x)=g_{0}x+g_{1}\in \Gamma.
\end{eqnarray}
Furthermore, $W$ is a quasi $\<\Gamma\cdot S\>$-module-at-infinity
with $\Phi$ being the identity map. \et

\begin{proof} For $g(x),h(x)\in \Gamma$, we have
$$g(x_{1})-h(x_{2})
=g_{0}\left(x_{1}-(g_{0}^{-1}h_{0}x_{2}+g_{0}^{-1}(h_{1}-g_{1}))\right)
=g_{0}(x_{1}-(g^{-1}h)(x_{2})).$$ With this, it is clear that
$\Gamma \cdot S$ is $\Gamma$-local. By Theorem \ref{tquasi-local},
$\Gamma\cdot S$ generates a vertex algebra $\<\Gamma\cdot S\>$
inside $\E^{o}(W)$ and $W$ is a quasi module-at-infinity for
$\<\Gamma\cdot S\>$. It follows from Lemma \ref{lconnection} and
induction that $\<\Gamma \cdot S\>$ is $\Gamma$-stable. In view of
Lemma \ref{lconnection}, $\<\Gamma \cdot S\>$ equipped with the
action of $\Gamma$ and with the group homomorphism $\phi$ is a
$\Gamma$-vertex algebra. Furthermore, for $\alpha(x),\beta(x)\in
\Gamma\cdot S$, as $Y_{W}(\alpha(x),x_{1})=\alpha(x_{1})$ and
$Y_{W}(\beta(x),x_{2})=\beta(x_{2})$, $Y_{W}(\alpha(x),x_{1})$ and
$Y_{W}(\beta(x),x_{2})$ are $\Gamma$-local. For $g\in \Gamma,\;
\alpha(x)\in \<\Gamma\cdot S\>$, we have
$$Y_{W}(R_{g}\alpha(x),x_{0})=Y_{W}(\alpha(g(x)),x_{0})=\alpha(g(x_{0}))
=Y_{W}(\alpha(x),g(x_{0})).$$ It follows from Lemma
\ref{lva-gamma-quasi-module} that $W$ is a quasi module-at-infinity
for $\<\Gamma \cdot S\>$ viewed as a $\Gamma$-vertex algebra.
\end{proof}

We shall need the following result:

 \bl{ldecomposition} Let $W$ be a
vector space and let $a(x),b(x)\in \E^{o}(W)$. Suppose that
\begin{eqnarray}\label{eab-bracket}
-[a(x_{1}),b(x_{2})]=\sum_{i=1}^{k}\sum_{j=0}^{r}\Psi_{i,j}(x_{2})
\frac{1}{j!}\left(\frac{\partial}{\partial x_{2}}\right)^{j}
x_{1}^{-1}\delta\left(\frac{\beta_{i}x_{2}}{x_{1}}\right),
\end{eqnarray}
where $\beta_{1},\dots,\beta_{k}$ are distinct nonzero complex
numbers with $\beta_{1}=1$ and $\Psi_{i,j}(x)\in \E^{o}(W)$.
Then $(a(x),b(x))$ is quasi local and $a(x)_{n}b(x)=0$ for $n>r$
and $a(x)_{n}b(x)=\Psi_{1,n}(x)$
for $0\le n\le r$.
\el

\begin{proof} Set $p(x,z)=(x-\beta_{1}z)^{r+1}\cdots
  (x-\beta_{k}z)^{r+1}$ and $q(x,z)=(x-\beta_{2})^{r+1}\cdots (x-\beta_{k}z)^{r+1}$.
{}From (\ref{eab-bracket}) we have
$p(x_{1},x_{2})[a(x_{1}),b(x_{2})]=0$. Thus $(a(x),b(x))$ is quasi
local. Furthermore, we have
\begin{eqnarray*}
p(x+x_{0},x)Y_{\E^{o}}(a(x),x_{0})b(x)
=\left(p(x_{1},x)a(x_{1})b(x)\right)|_{x_{1}=x+x_{0}}.
\end{eqnarray*}
In view of Lemma \ref{lsimple-factx} we have
\begin{eqnarray*}
& &x_{0}^{-1}\delta\left(\frac{x_{1}-x}{x_{0}}\right)
p(x_{1},x)b(x)a(x_{1})-
x_{0}^{-1}\delta\left(\frac{x-x_{1}}{-x_{0}}\right)
p(x_{1},x)a(x_{1})b(x)\\
& &\ \ \ \ =x_{1}^{-1}\delta\left(\frac{x+x_{0}}{x_{1}}\right)
p(x_{1},x)Y_{\E^{o}}(a(x),x_{0})b(x).
\end{eqnarray*}
Multiplying both sides by $x_{0}^{-r-1}$ and using delta-function
substitution we get
\begin{eqnarray*}
& &x_{0}^{-1}\delta\left(\frac{x_{1}-x}{x_{0}}\right)
q(x_{1},x)b(x)a(x_{1})-
x_{0}^{-1}\delta\left(\frac{x-x_{1}}{-x_{0}}\right)
q(x_{1},x)a(x_{1})b(x)\\
& &\ \ \ \ =x_{1}^{-1}\delta\left(\frac{x+x_{0}}{x_{1}}\right)
q(x_{1},x)Y_{\E^{o}}(a(x),x_{0})b(x).
\end{eqnarray*}
Taking $\Res_{x_{0}}$ we get
\begin{eqnarray}\label{efirst}
-q(x_{1},x)[a(x_{1}),b(x)]&=&
\Res_{x_{0}}x_{1}^{-1}\delta\left(\frac{x+x_{0}}{x_{1}}\right)
q(x_{1},x)Y_{\E^{o}}(a(x),x_{0})b(x)\nonumber\\
&=&\sum_{j\ge 0}q(x_{1},x)a(x)_{j}b(x)
\frac{1}{j!}\left(\frac{\partial}{\partial x}\right)^{j}
x_{1}^{-1}\delta\left(\frac{x}{x_{1}}\right).
\end{eqnarray}
On the other hand, {}from (\ref{eab-bracket}) we have
\begin{eqnarray}\label{esecond}
-q(x_{1},x)[a(x_{1}),b(x)]=\sum_{j=0}^{r}q(x_{1},x)\Psi_{1,j}(x)
\frac{1}{j!}\left(\frac{\partial}{\partial x}\right)^{j}
x_{1}^{-1}\delta\left(\frac{x}{x_{1}}\right).
\end{eqnarray}
 Assume
$$\sum_{j=0}^{s}q(x_{1},x)A_{j}(x)
\frac{1}{j!}\left(\frac{\partial}{\partial x}\right)^{j}
x_{1}^{-1}\delta\left(\frac{x}{x_{1}}\right)=0,$$ where $A_{j}(x)\in
\E^{o}(W)$ for $0\le j\le s$. As $q(x,1)$ and $(x-1)$ are relatively
prime, we have $1=q(x,1)f(x)+(x-1)^{s+1}g(x)$ for some $f(x),g(x)\in
\C[x]$. Then
\begin{eqnarray*}
& &\sum_{j=0}^{s}A_{j}(x) \frac{1}{j!}\left(\frac{\partial}{\partial
x}\right)^{j}
x_{1}^{-1}\delta\left(\frac{x}{x_{1}}\right)\\
&=&\sum_{j=0}^{s}A_{j}(x)\left(q(x_{1}/x,1)f(x_{1}/x)+(x_{1}/x-1)^{s+1}g(x_{1}/x)\right)
\frac{1}{j!}\left(\frac{\partial}{\partial x}\right)^{j}
x_{1}^{-1}\delta\left(\frac{x}{x_{1}}\right)
\\  &=&0,
\end{eqnarray*}
 which implies $A_{j}(x)=0$ for $0\le j\le s$.
Using this fact,  combining (\ref{efirst}) with (\ref{esecond}) we
obtain the desired relations.
\end{proof}

Next, we study certain infinite-dimensional Lie algebras including
the Lie algebra of pseudo-differential operators on the circle.
Let $\g$ be a (possibly infinite-dimensional) Lie algebra equipped
with a nondegenerate symmetric invariant bilinear form $\<\cdot,\cdot\>$.
Associated with the pair $(\g,\<\cdot,\cdot\>)$, one has an
(untwisted) affine Lie algebra
$$\hat{\g}=\g\otimes \C[t,t^{-1}]\oplus \C {\bf k},$$
where ${\bf k}$ is central and
$$[a\otimes t^{m},b\otimes
t^{n}]=[a,b]\otimes t^{m+n}+m\delta_{m+n,0}\<a,b\>{\bf k}$$
for
$a,b\in \g,\; m,n\in \Z$. Defining $\deg (\g\otimes t^{m})=-m$ for
$m\in \Z$ and $\deg {\bf k}=0$ makes $\hat{\g}$ a $\Z$-graded Lie
algebra. For $a\in \g$, form the generating function
 $$a(x)=\sum_{n\in \Z}(a\otimes t^{n})x^{-n-1}.$$

We say that a $\hat{\g}$-module $W$ is of {\em level} $\ell\in \C$
if ${\bf k}$ acts on $W$ as scalar $\ell$. A {\em vacuum vector} in
a $\hat{\g}$-module is a nonzero vector $v$ such that $(\g\otimes
\C[t])v=0$ and a {\em vacuum $\hat{\g}$-module} is a
$\hat{\g}$-module $W$ equipped with a vacuum vector which generates
$W$. Let $\ell$ be any complex number. Denote by $\C_{\ell}$ the
$1$-dimensional $(\g\otimes \C[t]\oplus \C {\bf k})$-module with
$\g\otimes \C[t]$ acting trivially and with ${\bf k}$ acting as
scalar $\ell$. Form the induced $\hat{\g}$-module
$$V_{\hat{\g}}(\ell,0)=U(\hat{\g})\otimes_{U(\g\otimes \C[t]
\oplus \C {\bf k})} \C_{\ell}.$$
Set ${\bf 1}=1\otimes 1$, which is a vacuum vector. {}From
definition, $V_{\hat{\g}}(\ell,0)$ is a universal vacuum
$\hat{\g}$-module of level $\ell$. Identify $\g$ as a subspace of
$V_{\hat{\g}}(\ell,0)$ through the linear map $a\rightarrow
a(-1){\bf 1}$. It is now well known (cf. \cite{fz}) that there
exists one and only one vertex-algebra structure on $V_{\hat{\g}}(\ell,0)$
with ${\bf 1}$ as the vacuum vector and with $Y(a,x)=a(x)$ for $a\in
\g$. Defining $\deg {\bf 1}=0$ makes $V_{\hat{\g}}(\ell,0)$ a
$\Z$-graded $\hat{\g}$-module and
the vertex algebra $V_{\hat{\g}}(\ell,0)$ equipped with this
$\Z$-grading is a $\Z$-graded vertex algebra.

Let $\Gamma$ be a subgroup of $\Aut(\g)$,
preserving the bilinear form $\<\cdot,\cdot\>$.
Each $g\in \Gamma$ canonically lifts to an automorphism of
the $\Z$-graded Lie algebra $\hat{\g}$. Then
$\Gamma$ acts on the vertex algebra $V_{\hat{\g}}(\ell,0)$
by automorphisms preserving the $\Z$-grading.
Let $\phi: \Gamma\rightarrow \C^{\times}$
be a group homomorphism. For $g\in \Gamma$, set
$$R_{g}=\phi(g)^{-L(0)}g\in GL(V_{\hat{\g}}(\ell,0)),$$
where $L(0)$ denotes the $\Z$-grading operator. This defines a
$\Gamma$-vertex-algebra structure on $V_{\hat{\g}}(\ell,0)$.

Consider the following completion of the $\Z$-graded affine Lie algebra
$\hat{\g}$:
\begin{eqnarray}
\hat{\g}(\infty)=\g\otimes \C((t^{-1}))\oplus \C {\bf k},
\end{eqnarray}
where
\begin{eqnarray}
[a\otimes p(t),b\otimes q(t)]=[a,b]\otimes
p(t)q(t)+\Res_{t}p'(t)q(t)\<a,b\>{\bf k}
\end{eqnarray}
for $a,b\in \g,\; p(t),q(t)\in \C((t^{-1}))$.

\bp{ptwisted-affine-comp}
Let $\g$ be a Lie algebra equipped with
a nondegenerate symmetric invariant bilinear form
 $\<\cdot,\cdot\>$ and let $\Gamma$ be a group acting on $\g$ by automorphisms
 preserving the bilinear form $\<\cdot,\cdot\>$ and satisfying the condition
 that for any $u,v\in \g$,
$$[gu,v]=0\;\;\mbox{ and }\;\; \<gu,v\>=0\;\;\;\mbox{ for all but
finitely many }g\in \Gamma.$$ Let $\Phi: \Gamma\rightarrow G$ be a
group homomorphism. Define a new bilinear multiplicative operation
$[\cdot,\cdot]_{\Gamma}$ on the vector space
$\hat{\g}(\infty)=\g\otimes \C((t^{-1}))\oplus \C {\bf k}$ by
\begin{eqnarray*}
& &[a\otimes p(t),{\bf k}]_{\Gamma}=0=[{\bf k},a\otimes p(t)]_{\Gamma}, \\
& &[a\otimes p(t),b\otimes q(t)]_{\Gamma}=\sum_{g\in
\Gamma}[ga,b]\otimes p(g(t))q(t)+
\Res_{t}p'(g(t))g'(t)q(t)\<ga,b\>{\bf k}\ \ \ \ \ \
\end{eqnarray*}
for $a,b\in \g,\; p(t),q(t)\in \C((t^{-1}))$. Then the subspace,
linearly spanned by the elements $$ga\otimes p(g(t))-a\otimes p(t)$$
for $g\in \Gamma,\;a\in \g,\;p(t)\in\C((t^{-1}))$, is a two-sided
ideal of the nonassociative algebra and the quotient algebra, which
we denote by $\hat{\g}(\infty)[\Gamma]$, is a Lie algebra.
\ep

\begin{proof} Let $\Gamma$ act on the Lie algebra $\hat{\g}(\infty)$ by
$$g(a\otimes p(t)+\lambda {\bf k})
=ga\otimes p(g(t))+\lambda {\bf k}\;\;\;\mbox{ for }g\in \Gamma,\;
a\in \g,\; p(t)\in \C((t^{-1})),\; \lambda\in \C.$$ It is
straightforward to see that $\Gamma$ acts on $\hat{\g}(\infty)$ by
automorphisms. We have
\begin{eqnarray*}
\sum_{g\in \Gamma}[g(a\otimes p(t)),b\otimes q(t)]
&=&\sum_{g\in \Gamma}[ga\otimes p(g(t)),b\otimes q(t)]\\
&=&\sum_{g\in \Gamma}[ga,b]\otimes p(g(t))q(t)
+\Res_{t}q(t)p'(g(t))g'(t)\<ga,b\>{\bf k}
\end{eqnarray*}
for $a,b\in \g,\; p(t),q(t)\in \C((t^{-1}))$ and
$$g(a\otimes p(t)+\lambda {\bf k})-(a\otimes p(t)+\lambda {\bf
k})=ga\otimes p(g(t))-a\otimes p(t).$$ Now the assertions follow
immediately from (\cite{li-gamodule}, Lemma 4.1).
\end{proof}

Let
$$\pi:\hat{\g}(\infty)\rightarrow \hat{\g}(\infty)[\Gamma] $$
be the natural map.
For $a\in \g$,  set
$$a_{\Gamma}(x)=\sum_{n\in \Z}\pi (a\otimes t^{n}) x^{-n-1}
\in \left(\hat{\g}(\infty)[\Gamma]\right)[[x,x^{-1}]].$$

\bl{lgenerating-relation}
For $g\in \Gamma,\; a\in \g$, we have
\begin{eqnarray}\label{e5.15}
(ga)_{\Gamma}(x)=\phi(g)a_{\Gamma}(g(x)).
\end{eqnarray}
For $a,b\in \g$, we have
\begin{eqnarray}\label{ecomm-liealgebra}
[a_{\Gamma}(x_{1}),b_{\Gamma}(x_{2})]=\sum_{g\in
\Gamma}[ga,b]_{\Gamma}(x_{2})
x_{1}^{-1}\delta\left(\frac{g(x_{2})}{x_{1}}\right) +\<ga,b\>{\bf
k}\frac{\partial}{\partial
x_{2}}x_{1}^{-1}\delta\left(\frac{g(x_{2})}{x_{1}}\right),
\end{eqnarray}
where $g(x)=\Phi(g)(x)\in G$. \el

\begin{proof}
Let $g(x)=g_{0} x+g_{1}\in G$, where $g_{0}\in \C^{\times},\;
g_{1}\in \C$. Then $g^{-1}(x)=g_{0}^{-1}(x-g_{1})$ and
\begin{eqnarray}
x^{-1}\delta\left(\frac{g^{-1}(t)}{x}\right)
=x^{-1}\delta\left(\frac{t-g_{1}}{g_{0}x}\right)
=g_{0} t^{-1}\delta\left(\frac{g_{0} x+g_{1}}{t}\right)
=g_{0} t^{-1}\delta\left(\frac{g(x)}{t}\right).
\end{eqnarray}
{}From  definition we have
$\pi(ga\otimes t^{n})=\pi(a\otimes (g^{-1}(t))^{n})$
for $g\in \Gamma,\;a\in \g,\; n\in \Z$. Then we get
\begin{eqnarray*}
(ga)_{\Gamma}(x)&=&\sum_{n\in\Z}\pi(ga\otimes t^{n})x^{-n-1}
= \sum_{n\in\Z}\pi(a\otimes
(g^{-1}(t))^{n})x^{-n-1}\\
&=&\pi\left(a\otimes x^{-1}\delta\left(\frac{g^{-1}(t)}{x}\right)\right)\\
&=&g_{0} \pi\left(a\otimes
t^{-1}\delta\left(\frac{g(x)}{t}\right)\right)\\
&=&g_{0}a_{\Gamma}(g(x)).
\end{eqnarray*}
proving (\ref{e5.15}). As for (\ref{ecomm-liealgebra}), notice that
$$\sum_{m,n\in \Z}g(t)^{m}t^{n}x_{1}^{-m-1}x_{2}^{-n-1}
=x_{1}^{-1}\delta\left(\frac{g(t)}{x_{1}}\right)
x_{2}^{-1}\delta\left(\frac{t}{x_{2}}\right)
=x_{1}^{-1}\delta\left(\frac{g(x_{2})}{x_{1}}\right)
x_{2}^{-1}\delta\left(\frac{t}{x_{2}}\right) $$
and
\begin{eqnarray*}
\Res_{t}\sum_{m,n\in \Z}mg(t)^{m-1}g'(t)t^{n}x_{1}^{-m-1}x_{2}^{-n-1}
&=&-\Res_{t}\sum_{m,n\in \Z}ng(t)^{m}t^{n-1}x_{1}^{-m-1}x_{2}^{-n-1}\\
&=&\Res_{t}\frac{\partial}{\partial x_{2}}
x_{1}^{-1}\delta\left(\frac{g(t)}{x_{1}}\right)
t^{-1}\delta\left(\frac{x_{2}}{t}\right)\\
&=&\frac{\partial}{\partial x_{2}}
x_{1}^{-1}\delta\left(\frac{g(x_{2})}{x_{1}}\right).
\end{eqnarray*}
Now (\ref{ecomm-liealgebra}) follows.
\end{proof}

Now we are in a position to present our main result of this section,
which is an analogue of a theorem of \cite{li-gamodule}:

\bt{tquasi-module-lie} Let $\g,\<\cdot,\cdot\>, \Gamma, \Phi$ be
given as in Proposition \ref{ptwisted-affine-comp} and let $W$ be
any $\hat{\g}(\infty)[\Gamma]$-module of level $\ell\in \C$ such
that $a_{\Gamma}(x)\in \E^{o}(W)$ for $a\in \g$. Then on $W$ there
exists one and only one structure of a quasi module-at-infinity for
$V_{\hat{\g^{o}}}(-\ell,0)$ viewed as a $\Gamma$-vertex algebra with
$Y_{W}(a,x)=a_{\Gamma}(x)$ for $a\in \g$, where $\g^{o}$ denotes the
opposite Lie algebra of $\g$. \et

\begin{proof} Set
$$U=\{ a_{\Gamma}(x)\;|\; a\in \g\}\subset \E^{o}(W).$$
For $a,b\in \g$, let $g_{1},\dots,g_{r}\in \Gamma$ such that
$[ga,b]=0$ and $\<ga,b\>=0$ for $g\notin \{g_{1},\dots,g_{r}\}$.
It follows from (\ref{ecomm-liealgebra}) that
$$(x_{1}-g_{1}(x_{2}))^{2}\cdots
(x_{1}-g_{r}(x_{2}))^{2}[a_{\Gamma}(x_{1}),b_{\Gamma}(x_{2})]=0.$$
Thus $U$ is a $\Gamma$-local subspace of $\E^{o}(W)$. {}From
(\ref{e5.15}), $\Gamma\cdot U=U$. By Theorem \ref{tmain-2}, $U$
generates a $\Gamma$-vertex algebra $\<U\>$ with $W$ as a quasi
module-at-infinity with $Y_{W}(\alpha(x),x_{0})=\alpha(x_{0})$.
Combining (\ref{ecomm-liealgebra}) with Lemma \ref{ldecomposition}
we get
\begin{eqnarray}
a_{\Gamma}(x)_{0}b_{\Gamma}(x)=-[a,b]_{\Gamma}(x),\;\;
a_{\Gamma}(x)_{1}b_{\Gamma}(x)=-\ell \<a,b\>1_{W}, \;\mbox{ and
}\;a_{\Gamma}(x)_{n}b_{\Gamma}(x)=0
\end{eqnarray}
for $n\ge 2$. In view of the universal property (cf. \cite{pr},
\cite{li-gamodule}) of $V_{\hat{\g^{o}}}(-\ell,0)$, there exists a
(unique) vertex-algebra homomorphism from
$V_{\hat{\g^{o}}}(-\ell,0)$ to $\<U\>$, sending $a$ to
$a_{\Gamma}(x)$ for $a\in \g$. Consequently, $W$ is a quasi
module-at-infinity for $V_{\hat{\g^{o}}}(-\ell,0)$ viewed as a
vertex algebra. Furthermore, for $g\in \Gamma,\; a\in \g$ we have
$$Y_{W}(R_{g}a,x)=\phi(g)^{-1}(ga)_{\Gamma}(x)
=a_{\Gamma}(g(x))=Y_{W}(a,g(x)).$$ As $\g$ generates
$V_{\widehat{\g^{o}}}(-\ell,0)$ as a vertex algebra, it follows from
Lemma \ref{lva-gamma-quasi-module} that $W$ is a quasi
module-at-infinity for $V_{\hat{\g^{o}}}(-\ell,0)$ viewed as a
$\Gamma$-vertex algebra.
\end{proof}

Now, let $\Gamma$ be any abstract group. As in \cite{li-gamodule},
we define an associative algebra $gl_{\Gamma}$ with a $\C$-basis
consisting of symbols $E_{\alpha,\beta}$ for $\alpha,\beta\in
\Gamma$ and with
$$E_{\alpha,\beta}E_{\mu,\nu}=\delta_{\beta,\mu}E_{\alpha,\nu}
\ \ \ \mbox{  for }\alpha,\beta,\mu,\nu\in \Gamma,$$ and we equip
$gl_{\Gamma}$ with a nondegenerate symmetric associative bilinear
form defined by
$$\<E_{\alpha,\beta},E_{\mu,\nu}\>=\delta_{\alpha,\nu}\delta_{\beta,\mu}
\;\;\;\mbox{ for }\alpha,\beta,\mu,\nu\in \Gamma.$$
Defining $T_{\alpha}\in
GL(gl_{\Gamma})$ for $\alpha\in \Gamma$ by
$$T_{\alpha}E_{\mu,\nu}=E_{\alpha \mu,\alpha\nu}
\;\;\;\mbox{ for }\alpha,\beta,\mu,\nu\in \Gamma$$ (cf. \cite{gkk})
we have a group action of $\Gamma$ on $gl_{\Gamma}$  by
automorphisms preserving the bilinear form. We can also view
$gl_{\Gamma}$ as a Lie algebra with $\<\cdot,\cdot\>$ an invariant
bilinear form. Furthermore, for any $\alpha,\beta,\mu,\nu\in
\Gamma$, we have
$$[T_{g}E_{\alpha,\beta},E_{\mu,\nu}]=0\;\;\;\mbox{ and }
\;\;\; \< T_{g}E_{\alpha,\beta},E_{\mu,\nu}\>=0$$ for all but
finitely many $g\in \Gamma$. Associated with the pair
$(gl_{\Gamma},\<\cdot,\cdot\>)$, we have an (untwisted) affine Lie
algebra $\widehat{gl_{\Gamma}}$ and its completion
$\widehat{gl_{\Gamma}}(\infty)$.

Let $\Phi: \Gamma\rightarrow G$ be a group homomorphism. For
$\alpha\in \Gamma$, we set
$$\alpha (x)=\Phi(\alpha)=\alpha_{0}x+\alpha_{1}\in G,$$
where $\alpha_{0},\alpha_{1}\in \C$ with $\alpha_{0}\ne 0$.
{}From Proposition \ref{ptwisted-affine-comp}
we have a Lie algebra $\widehat{gl_{\Gamma}}(\infty)[\Gamma]$.

For $\alpha,\beta\in \Gamma$, we have
\begin{eqnarray}\label{egeneral-comm}
& &[E_{\alpha,e}(x_{1}),E_{\beta,e}(x_{2})]_{\Gamma}\nonumber\\
&=&\sum_{g\in
\Gamma}[E_{g\alpha,g},E_{\beta,e}](x_{2})
x_{1}^{-1}\delta\left(\frac{g(x_{2})}{x_{1}}\right)
+\<E_{g\alpha,g},E_{\beta,e}\>{\bf k}\frac{\partial}{\partial
x_{2}}x_{1}^{-1}\delta\left(\frac{g(x_{2})}{x_{1}}\right)\nonumber\\
&=&E_{\beta\alpha,e}(x_{2})
x_{1}^{-1}\delta\left(\frac{\beta(x_{2})}{x_{1}}\right)
-E_{\beta,\alpha^{-1}}(x_{2})
x_{1}^{-1}\delta\left(\frac{\alpha^{-1}(x_{2})}{x_{1}}\right)
\nonumber\\
& &\ \ \ \ +\delta_{\alpha\beta,e}{\bf k} \frac{\partial}{\partial
x_{2}}x_{1}^{-1}\delta\left(\frac{\beta(x_{2})}{x_{1}}\right).
\end{eqnarray}
For $\alpha\in \Gamma$, denote by $A_{\alpha}(x)$ the image of
$E_{\alpha,e}(x)$ in
$\left(\widehat{gl_{\Gamma}}(\infty)[\Gamma]\right)[[x,x^{-1}]]$.
Note that
$$E_{\beta,\alpha^{-1}}(x)=(T_{\alpha^{-1}}E_{\alpha\beta,e})(x)
=\alpha_{0}E_{\alpha\beta,e}(\alpha^{-1}(x)). $$ By
(\ref{egeneral-comm}) we have
\begin{eqnarray}
[A_{\alpha}(x_{1}),A_{\beta}(x_{2})]
&=&A_{\beta\alpha}(x_{2})
x_{1}^{-1}\delta\left(\frac{\beta(x_{2})}{x_{1}}\right)
-\alpha_{0}A_{\alpha\beta}(\alpha^{-1}(x_{2}))
x_{1}^{-1}\delta\left(\frac{\alpha^{-1}(x_{2})}{x_{1}}\right)
\nonumber\\
& &\ \ \ \ \ \ +\delta_{\alpha\beta,e}{\bf k} \frac{\partial}{\partial
x_{2}}x_{1}^{-1}\delta\left(\frac{\beta(x_{2})}{x_{1}}\right).
\end{eqnarray}
Let $W$ be any $\widehat{gl_{\Gamma}}(\infty)[\Gamma]$-module of
level $\ell\in \C$ such that $E_{\alpha,\beta}(x)\in \E^{o}(W)$ for
all $\alpha,\beta\in \Gamma$. In view of Theorem
\ref{tquasi-module-lie}, there exists one and only one quasi
module-at-infinity structure on $W$ for the $\Gamma$-vertex algebra
$V_{\widehat{gl_{\Gamma}^{o}}}(-\ell,0)$.

\bex{example-pseudo} {\em Consider the special case with
$\Gamma=\Z$. (The associative algebra $gl_{\Z}$ is simply
$gl_{\infty}$, the associative algebra of doubly infinite complex
matrices with only finite many nonzero entries.) Let $\phi$ be the
group embedding of $\Z$ into $G$ defined by $\phi(n)(x)=x+n$ for
$n\in \Z$. We have
\begin{eqnarray}
[A_{m}(x_{1}),A_{n}(x_{2})]
&=&A_{m+n}(x_{2})x_{1}^{-1}\delta\left(\frac{x_{2}+n}{x_{1}}\right)
-A_{m+n}(x_{2}-m)x_{1}^{-1}\delta\left(\frac{x_{2}-m}{x_{1}}\right)
\nonumber\\
& &\hspace{1cm}+\delta_{m+n,0}{\bf k} \frac{\partial}{\partial
x_{2}}x_{1}^{-1}\delta\left(\frac{x_{2}+n}{x_{1}}\right)\nonumber\\
&=&A_{m+n}(x_{2})x_{1}^{-1}\delta\left(\frac{x_{2}+n}{x_{1}}\right)
-A_{m+n}(x_{1})x_{2}^{-1}\delta\left(\frac{x_{1}+m}{x_{2}}\right)
\nonumber\\
& &\hspace{1cm}+\delta_{m+n,0}{\bf k} \frac{\partial}{\partial
x_{2}}x_{1}^{-1}\delta\left(\frac{x_{2}+n}{x_{1}}\right)
\end{eqnarray}
for $m,n\in \Z$. This is exactly the relation of the Lie algebra of
pseudo-differential operators on the circle. This Lie algebra
(without central extension) was studied in \cite{gkk} in terms of a
notion called $\Gamma$-conformal algebra. Note that this Lie algebra
does not admit nontrivial modules of highest weight type. } \eex

\end{document}